\documentclass[leqno]{article}
\usepackage{verbatim, amsmath, amscd, amssymb}
\usepackage[all]{xy}
\newcommand{\cB}{{\mathcal B}}

\newcommand{\cD}{{\mathcal D}}
\newcommand{\cE}{{\mathcal E}}
\newcommand{\cF}{{\mathcal F}}

\newcommand{\cL}{{\mathcal L}}
\newcommand{\cM}{{\mathcal M}}
\newcommand{\cMod}{{\mathcal M}\!{\it od}}

\newcommand{\cO}{{\mathcal O}}
\newcommand{\cP}{{\mathcal P}}

\newcommand{\frg}{\mathfrak g}

\newcommand{\frp}{\mathfrak p}

\newcommand{\rd}{\mathrm{d}}

\newcommand{\rD}{\mathrm{D}}

\newcommand{\rG}{\mathrm{G}}\newcommand{\rH}{\mathrm{H}}

\newcommand{\rN}{\mathrm{N}}
\newcommand{\rR}{\mathrm{R}}

\newcommand{\rSO}{\mathrm{SO}}

\newcommand{\bbB}{\mathbb B}

\newcommand{\bbF}{\mathbb F}

\newcommand{\bbP}{\mathbb P}
\newcommand{\bbN}{\mathbb N}
\newcommand{\bbQ}{\mathbb Q}

\newcommand{\bbZ}{\mathbb Z}

\newcommand{\bfR}{\mathbf R}

\newcommand{\Ad}{\mathrm{Ad}}
\newcommand{\ad}{\operatorname{ad}}

\newcommand{\coh}{\mathrm{coh}}

\newcommand{\ch}{\mathrm{ch}\,}

\newcommand{\diag}{\mathrm{diag}}

\newcommand{\cDiff}{\mathcal{D}if\!f}

\newcommand{\Ext}{\mathrm{Ext}}

\newcommand{\GL}{\mathrm{GL}}
\newcommand{\Gr}{\mathrm{Gr}}

\newcommand{\id}{\mathrm{id}}

\newcommand{\ind}{\mathrm{ind}}

\newcommand{\rad}{\mathrm{rad}}

\newcommand{\SL}{\mathrm{SL}}
\newcommand{\soc}{\mathrm{soc}}

\newcommand{\wt}{\mathrm{wt}}

\newcommand{\Mod}{\mathbf{Mod}}

\newcommand{\lbr}{\begin{bmatrix}}
\newcommand{\rbr}{\end{bmatrix}}
\newcommand{\for}{\bigcirc\kern-2.6ex \because}
\newcommand{\forb}{\bigcirc\kern-2.8ex \because}
\newcommand{\forbb}{\bigcirc\kern-3.0ex \because}
\newcommand{\forbbb}{\bigcirc\kern-3.1ex \because}

\newcommand{\tto}{\twoheadrightarrow}

\newcommand\pf{\noindent {\bf Proof:  }}
\newtheorem{thm}{Theorem:}

\newtheorem{prop}{Proposition:}

\newtheorem{cor}{Corollary:}

\begin{document}
\large
\title{
{\bf 
On the Frobenius direct image of the structure sheaf of a homogeneous projective variety 
}
\thanks{supported in part by JSPS Grants in Aid for Scientific Research 15K04789 
}
\author{
K\textsc{aneda} Masaharu
\\
Osaka City University
\\
Department of Mathematics
\\
kaneda@sci.osaka-cu.ac.jp
}
}
\maketitle

\begin{abstract}
We present an example of 
a homogeneous projective variety 
the Frobenius direct image of the structure sheaf of which 
has nonvanishing self extension.

\end{abstract}

On some homogeneous projective spaces
$G/P$
in large positive
characteristic
we have found, for the projective spaces
\cite{09}, the quadrics
\cite{14}, Grassmannians
$\Gr(2,n)$
\cite{17}, 
and when $G$ is the special linear group  of degree 3
\cite{HKR}, the symplectic group of degree 4
\cite{AK00}, \cite{07}, or 
when
$G$ is in type $\rG_2$ and $P$ a maximal parabolic subgroup
\cite{KY},
a Karoubian complete
strongly exceptional
collection 
of coherent modules 
$\cE(w)$
over $G/P$,
parametrized by 
the minimal length representatives 
$w$
of the cosets of the Weyl group of $G$ by that of $P$,
as subquotients of the Frobenius direct image of the structure sheaf of $G/P$,
which are all defined over $\bbZ$.
Except for the cases of $\Gr(2,n)$, $n\geq4$,
and when $G$ is in type $G_2$ with $P$ associated to a short simple root,
we also know that
the Frobenius direct image of the structure sheaf is a direct sum of copies of those
$\cE(w)$'s;
the case for the quadrics is due to
Langer
\cite{La}.
In this paper we determine an extra summand
in that exceptional case in type $\rG_2$,
and find that the summand causes nontrivial self
extension of the Frobenius direct image of the structure sheaf.
There follows 
nonvanishing of the 1st cohomology of the sheaf 
of rings of 
small differential operators on
$G/P$ in this setting.
The sheaf 
of rings of 
small differential operators is
the first term of the $p$-filtration 
\cite{Haa}
of 
the sheaf $\cDiff_{G/P}$
of rings of differential operators \cite{EGA}, 
and is a central reducion of the sheaf $\cD_{G/P}^{(0)}$
of rings of arithmetic differential operators of level 0
\cite{B} 
which is called the sheaf of rings of cristaline 
differential operators in \cite{BMR}.
One may recall that
Kashiwara and Lauritzen 
\cite{KaLa}
found
the nonvanishing of higher cohomolgy of $\cDiff_{G/P}$ for  
$G/P=\Gr(2,5)$,
while that
the vanishing of the higher cohomology of
$\cD_{G/P}^{(0)}$
in general
holds thanks to
Bezrukavnikov, Mirkovic, and Rumynin \cite{BMR}.

In more details,
let
$G$ be a simple algebraic group over an algebraically closed field of characteristic
at least 11, and $P$ the standard parabolic subgroup of $G$ associated to a short simple root $\alpha_1$.
Let $W$ be the Weyl group of $G$ 
with simple reflections $s_1$ and $s_2$,
and let 
$W^P=\{w\in W|w\alpha_1>0)=
\{
e,s_2,s_1s_2,s_2s_1s_2,s_1s_2s_1s_2,s_2s_1s_2s_1s_2\}$.
For each $w\in W^P$ let
$L(w)$ be the simple $G_1$-module of highest weight $w\rho-\rho$, $\rho$ a half sum of the positive roots
and $G_1$ the Frobenius kernel of $G$.
Let $\frp$ be the Lie algebra of $P$ under the adjoint represention, and let
$\cL_{G/P}(\frp)(-1)$ be the sheaf over $G/P$ associated to $\frp$ with Serre-twist $\cO_{G/P}(-1)$.

\begin{thm}
The Frobenius direct image 
$F_*\cO_{G/P}$
of the structure sheaf of $G/P$ decomposes into a direct sum of indecomposable sheaves
\begin{multline*}
\cE(e)\otimes
L(e)
\oplus
\cE(s_2)\otimes\{L(s_2)\oplus
L(s_1s_2s_1s_2)\}
\oplus
\cE(s_1s_2)\otimes
L(s_1s_2)\oplus
\cL_{G/P}(\frp)(-1)\otimes
L(e)
\oplus
\\
\cE(s_2s_1s_2)\otimes
L(s_2s_1s_2)
\oplus
\cE(s_1s_2s_1s_2)\otimes
L(s_1s_2s_1s_2)
\oplus
\cE(s_2s_1s_2s_1s_2)\otimes
\{L(s_2s_1s_2s_1s_2)\oplus
L(s_1s_2)\}.
\end{multline*}
The $\cE(w)$, $w\in W^P$, are all locally free sheaves of finite rank,
defined over $\bbZ$, and form a Karoubian complete strongly exceptional collection 
in the bounded derived category of coherent sheaves on $G/P$
such that
$\Mod_\cP(\cE(x),\cE(y))\ne0$ iff $x\geq y$ in the Chevalley-Bruhat order.
However,
$\Ext^1_\cP(F_*\cO_{G/P},F_*\cO_{G/P})\ne0$.

\end{thm}


\setcounter{equation}{0}
\begin{center}
$1^\circ$
{\bf 
Structure of the $G_1P$-Verma module
}

\end{center}

\setcounter{equation}{0}
\noindent
(1.1)
Let $\Bbbk$ 
be an algebraically closed field of positive
characteristic
$p$, $G$ a simple algebraic group over $\Bbbk$ in type $\rG_2$,
$B$ a Borel subgroup of $G$, $T$ a maximal torus of $B$,
$R$ the root system of $G$ relative to $T$,
$R^+$
the positive system of $R$ such that the roots of $B$ are $-R^+$,
and 
$R^s=\{\alpha_1,\alpha_2\}$ the set of simple roots with $\alpha_1$ short.
Let
$\Lambda$ be the character group of $T$, $\Lambda^+$ the set of dominant weights with the fundamental weights $\varpi_1$ and
$\varpi_2$;
$\langle\varpi_i,\alpha_j^\vee\rangle=\delta_{ij}$
$\forall i,j\in[1,2]$
with simple coroots
$\alpha_j^\vee$.
We partially order $\Lambda$ by 
$R^+$ such that
$\lambda\geq\mu$ iff $\lambda-\mu\in\sum_{\alpha\in R^+}\bbN\alpha$.
Let
$W$ be the Weyl group of $G$ with the simple reflections
$s_i$ associated to the simple root
$\alpha_i$, $i\in[1,2]$.
Let $P$ denote the standard parabolic subgroup of $G$ associated to the short simple root
$\alpha_1$ with the Weyl group $W_P$.
Let
$W^P=\{w\in W|w\alpha_1>0\}$ the set of minimal length representatives of
$W/W_P$.
Let $G_1$ be the Frobenius kernel of $G$ and let $\hat\nabla_P=\ind_P^{G_1P}$ be the induction functor from the category of $P$-modules to the category of $G_1P$-modules.
Let $\hat\nabla_P(\varepsilon)$
be the $G_1P$-Verma module of highest weight 0 induced from the trivial 1-dimensional $P$-module $\varepsilon$.
For $\lambda\in\Lambda^+$
we let
$L(\lambda)$ denote the simple $G$-module of highest weight $\lambda$.
We write each
$\mu\in\Lambda$ as a sum
$\mu=\mu^0+p\mu^1$ with
$\langle\mu^0,\alpha_i^\vee\rangle\in[0,p[$, $i\in\{1, 2\}$.
For $w\in W$ and $\mu\in\Lambda$ we let
$w\bullet\mu=w(\mu+\rho)-\rho$
with $\rho=\frac{1}{2}\sum_{\alpha\in R^+}\alpha=\varpi_1+\varpi_2$.
Put, in particular,
$L(w)=L((w\bullet0)^0)$,
which remains simple as $G_1$-module. 
For a $P$-module $M$ we let $M^{[1]}$ denote the Frobenius twist of $M$
\cite[II.3.16]{J}.
Unless otherwise specified, $\otimes$ will stand for the tensor product over $\Bbbk$.

We consider the geometric Frobenius morphism
$F:G/P\to G/P$ using the $\bbF_p$-form of $G/P$.
It factors through the natural morphism
$q:G/P\to G/G_1P$ to induce an isomorphism
$G_1P\to G/P$, so 
the Frobenius direct image $F_*\cO_{G/P}$ of the structure sheaf $\cO_{G/P}$ of $G/P$
may be identified with
the sheaf $\cL_{G/G_1P}(\hat\nabla_P(\varepsilon))$
over $G/G_1P$
associated to the $G_1P$-module
$\hat\nabla_P(\varepsilon)$.
Thus the structure of
$G_1P$-module on $\hat\nabla_P(\varepsilon)$ controls 
$G$-equivariantly
the structure of
$F_*\cO_{G/P}$.
Throughout the rest of the
paper we will assume 
$p\geq11$ so that Lusztig's conjecture for the irreducible characters for $G$ and $G_1T$ hold
\cite[D]{J}, which enables us to compute the $G_1T$-socle series of $\hat\nabla_P(\varepsilon)$ 
by the formula
\cite[5.2]{AbK}
using the periodic Kazhdan-Lusztig polynomials
\cite{L80}
and
\cite{Kato}.

\setcounter{equation}{0}
\noindent
(1.2)
Recall the $G_1T$-socle series of $\hat\nabla_P(\varepsilon)$:
$0=\soc^0\hat\nabla_P(\varepsilon)<\soc\hat\nabla_P(\varepsilon)=\soc^1\hat\nabla_P(\varepsilon)
<
\soc^2\hat\nabla_P(\varepsilon)
<\dots<
\soc^6\hat\nabla_P(\varepsilon)
=
\hat\nabla_P(\varepsilon)$.
Put
$\soc_i\hat\nabla_P(\varepsilon)=
(\soc^i\hat\nabla_P(\varepsilon))/
(\soc^{i-1}\hat\nabla_P(\varepsilon))$, $i\in[1,6]$.
One has a direct sum decomposition
$\soc_i\hat\nabla_P(\varepsilon)=
\coprod_{w\in W}L(w)\otimes
G_1\Mod(L(w),\soc_i\hat\nabla_P(\varepsilon))$.
As $G_1$ is normal in $G$, the decomposition holds as $G_1P$-modules.
Untwisting the Frobenius
or by the Frobenius contraction,
put
$\soc^1_{i,w}=G_1\Mod(L(w),\soc_i\hat\nabla_P(\varepsilon))^{[-1]}$
\cite[II.3.16]{J}/\cite{GK}.
Let $\nabla=\ind_B^G$
(resp. $\nabla^P=\ind_B^P$)
denote the induction functor from the category of $B$-modules to the category of
$G$-
(resp. $P$-) modules.
One has 
from \cite[4.7, 4.8]{KY}, see (A.1) in the appendix,
\begin{align*}
\soc_1\hat\nabla_P(\varepsilon)&\simeq
L(e)\otimes\varepsilon,
\\
\soc_2\hat\nabla_P(\varepsilon)&\simeq
L(s_2)\otimes
(-\varpi_2)^{[1]},
\\
\soc_3\hat\nabla_P(\varepsilon)&\simeq
L(s_1s_2)\otimes
\{(-\varpi_2)\otimes\ker(\nabla(\varpi_1)\twoheadrightarrow\nabla^P(\varpi_1))\}^{[1]}
\\
&\hspace{1cm}
\oplus
L(s_1s_2s_1s_2)\otimes
(-\varpi_2)^{[1]}
\oplus
L(e)\otimes(\soc_{3,e}^1)^{[1]},
\\
\soc_4\hat\nabla_P(\varepsilon)&\simeq
L(s_2s_1s_2)\otimes
\nabla^P(\varpi_1-2\varpi_2)^{[1]},
\\
\soc_5\hat\nabla_P(\varepsilon)&\simeq
L(s_1s_2s_1s_2)\otimes
\{(-2\varpi_2)\otimes(\nabla(\varpi_1)/\nabla^P(\varpi_1-\varpi_2))\}^{[1]}
\\
&\hspace{1cm}
\oplus
L(s_1s_2)\otimes(-2\varpi_2)^{[1]},
\\
\soc_6\hat\nabla_P(\varepsilon)&\simeq
L(w^P)\otimes
(-2\varpi_2)^{[1]},
\end{align*}
where
$\nabla(\varpi_1)\twoheadrightarrow\nabla^P(\varpi_1)$ is a unique epi of $P$-modules
and
$w^P=s_2s_1s_2s_1s_2$ is the longest element of
$W^P$.
We will determine the $P$-module
$\soc_{3,e}^1$ left open in
\cite{KY},
which will play the main role of the paper.

\setcounter{equation}{0}
\noindent
(1.3)
By the weight consideration
$\soc^1_{3,e}$ admits a $P$-module filtration
$0<-2\varpi_2=M^1<M^2<M^3<M^4=\soc^1_{3,e}$ such that $M^2/M^1\simeq\nabla^P(3\varpi_1-3\varpi_2)$,
$M^3/M^2\simeq-\varpi_2$,
and
$M^4/M^3\simeq\nabla^P(2\varpi_1-2\varpi_2)$.
We will denote a module with a filtration with subquotients $M_r,\dots, M_1$ from the top by
\begin{tabular}{|c|}
\hline
$M_r$
\\
\hline
\vdots
\\
\hline
$
M_1$
\\
\hline
\end{tabular}.
Thus, $\soc^1_{3,e}=
\begin{tabular}{|c|}
\hline
$\nabla^P(2\varpi_1-2\varpi_2)$
\\
\hline
$-\varpi_2$
\\
\hline
$\nabla^P(3\varpi_1-3\varpi_2)$
\\
\hline
$-2\varpi_2$
\\
\hline
\end{tabular}$.
Put $\soc^i=\soc^i\hat\nabla_P(\varepsilon)$,
$i\in[1,6]$.
Let
also $\hat\nabla=\ind_B^{G_1B}$ denote
the induction functor from the category of $B$-modules to the category of $G_1B$-modules.
We let $?^*$ denote the $\Bbbk$-linear dual.

Just suppose the extension 
$M^2$
of $\nabla^P(3\varpi_1-3\varpi_2)$ by 
$-2\varpi_2$ splits.
As $-3\varpi_1$ is the lowest weight of $\nabla^P(3\varpi_1-3\varpi_2)$ there would be a $G_1B$-submodule $M$
of $\soc^3$ containing $\soc^2$ such that
$M/\soc^2\simeq-3p\varpi_1$, and hence an exact sequence
$0\to \soc^2\otimes3p\varpi_1\to
M\otimes3p\varpi_1\to\varepsilon\to0$.
Applying the induction functor
$\ind_{G_1B}^G$ to the sequence
would then induce an exact sequence of $G$-modules
\begin{align}
\ind_{G_1B}^G(M\otimes3p\varpi_1)\to
\ind_{G_1B}^G(\varepsilon)\to
\rR^1\ind_{G_1B}^G(\soc^2\otimes3p\varpi_1)
\end{align}
with
\begin{align*}
\rR^1\ind_{G_1B}^G
&
(\soc^2\otimes3p\varpi_1)
=
\rR^1\ind_{G_1B}^G
(\begin{tabular}{|c|}
\hline
$L(s_2)\otimes(-\varpi_2)^{[1]}$
\\
\hline
$\varepsilon$
\\
\hline
\end{tabular}
\otimes3p\varpi_1)
\\
&\simeq
\rR^1\ind_{G_1B}^G
(\begin{tabular}{|c|}
\hline
$L(s_2)\otimes(3\varpi_1-\varpi_2)^{[1]}$
\\
\hline
$3p\varpi_1$
\\
\hline
\end{tabular}
).
\end{align*}
One has
\begin{align*}
\rR^1\ind_{G_1B}^G(3p\varpi_1)
&\simeq
\rR^1\ind_{B}^G(3\varpi_1)^{[1]} 
\quad\text{by \cite[II.3.19.3]{J}}
\\
&=0
\quad\text{by Kempf's vanishing
\cite[II.B.4]{J}},
\end{align*}
and also
\begin{align*}
\rR^1\ind_{G_1B}^G
&
(L(s_2)\otimes(3\varpi_1-\varpi_2)^{[1]})
\simeq
L(s_2)\otimes\rR^1\ind_{G_1B}^G((3\varpi_1-\varpi_2)^{[1]})
\\
&\hspace{5cm}
\text{by the tensor identity
\cite[I.4.8]{J}}
\\
&\simeq
L(s_2)\otimes\rR^1\ind_{G_1B}^G(3\varpi_1-\varpi_2)^{[1]}
=0
\quad\text{by Bott's theorem
\cite[II.5.5]{J}}.
\end{align*}
It follows that
$\rR^1\ind_{G_1B}^G
(\soc^2\otimes3p\varpi_1)
=0$.
On the other hand,
\begin{align*}
\ind_{G_1B}^G(M\otimes3p\varpi_1)
&\leq
\ind_{G_1B}^G(\hat\nabla_P(\varepsilon)\otimes3p\varpi_1)
\leq
\ind_{G_1B}^G(\hat\nabla(\varepsilon)\otimes3p\varpi_1)
\\
&\simeq
\ind_{G_1B}^G(\hat\nabla(3p\varpi_1))
\quad\text{by the tensor identity
again}
\\
&\simeq
\nabla(3p\varpi_1)
\quad\text{by the transitivity of inductions
\cite[I.3.5]{J}}.
\end{align*}
As $\ind_{G_1B}^G(\varepsilon)\simeq
\ind_B^G(\varepsilon)^{[1]}\simeq L(e)$, $L(e)$ would by (1) be a composition factor of $\nabla(3p\varpi_1)$, absurd
\cite[p. 150]{A86}.
Thus the extension is non-split.
Moreover,
\begin{align}
\Ext_P^1
&
(\nabla^P(3\varpi_1-3\varpi_2),
-2\varpi_2)
\simeq
\Ext_P^1
(\varepsilon,\nabla^P(3\varpi_1-3\varpi_2)^*\otimes(-2\varpi_2))
\\
\notag&
\simeq
\Ext_P^1
(\varepsilon,\nabla^P(-s_1(3\varpi_1-3\varpi_2))\otimes(-2\varpi_2))
\quad\text{as $\nabla^P(3\varpi_1-3\varpi_2)$ is $P$-simple}
\\
\notag&=
\Ext_P^1
(\varepsilon,\nabla^P(3\varpi_1)\otimes(-2\varpi_2))
\\
\notag&\simeq
\Ext_P^1
(\varepsilon,\nabla^P(3\varpi_1-2\varpi_2))
\quad\text{by the tensor identity}
\\
\notag&\simeq
\Ext_B^1
(\varepsilon,3\varpi_1-2\varpi_2)
\quad\text{by the Frobenius reciprocity
\cite[I.3.4]{J} and by Kempf}
\\
\notag&=
\Ext_B^1
(L(e),s_2\bullet0)
\\
\notag&\simeq
\Mod_G
(L(e),\rR^1\ind_B^G(s_2\bullet0))
\quad\text{by the Frobenius reciprocity and by Bott}
\\
\notag&\simeq
\Mod_G
(L(e),L(e))\simeq\Bbbk.
\end{align}
It follows that the non-split extension $M^2$
is unique up to isomorphism. 

Just suppose the extension 
$M^3$
of $-\varpi_2$ by $M^2$ is split.
There would then be a $G_1P$-submodule $M'$
of $\soc^3$ containing $\soc^2$ such that
$M'/\soc^2\simeq-\varpi_2$, and hence an exact sequence
$0\to \soc^2\otimes
p\varpi_2\to
M'\otimes
p\varpi_2\to\varepsilon\to0$.
Applying the induction functor
$\ind_{G_1P}^G$ to the sequence
would induce an exact sequence of $G$-modules
\[
\ind_{G_1P}^G(M'\otimes
p\varpi_2)\to
\ind_{G_1P}^G(\varepsilon)\to
\rR^1\ind_{G_1P}^G(\soc^2\otimes
p\varpi_2)
\]
with
\begin{align*}
\rR^1\ind_{G_1P}^G
&
(\soc^2\otimes
p\varpi_2)
=
\rR^1\ind_{G_1P}^G
(\begin{tabular}{|c|}
\hline
$L(s_2)\otimes(-\varpi_2)^{[1]}$
\\
\hline
$\varepsilon$
\\
\hline
\end{tabular}
\otimes
p\varpi_2)
\simeq
\rR^1\ind_{G_1P}^G
(\begin{tabular}{|c|}
\hline
$L(s_2)$
\\
\hline
$p\varpi_2$
\\
\hline
\end{tabular}
)
\\
&=0
\quad\text{by the tensor identity and by Kempf}.\end{align*}
Also,
$\ind_{G_1P}^G(M'\otimes
p\varpi_2)
\leq
\ind_{G_1P}^G(\hat\nabla_P(\varepsilon)\otimes
p\varpi_2)
\simeq
\ind_{G_1P}^G(\hat\nabla_P(p\varpi_2))
\simeq
\nabla(p\varpi_2)$.
Then
$L(e)$ would be a composition factor of $\nabla(p\varpi_2)$, absurd again
\cite{A86}.
Thus the extension $M^3$
is non-split.
Moreover,
there is a long exact sequence
\[
\dots\to
\Ext^1_P(-\varpi_2,
-2\varpi_2)\to
\Ext^1_P(-\varpi_2,
M^2)\to
\Ext^1_P(-\varpi_2,
\nabla^P(3\varpi_1-3\varpi_2))\to\dots
\]
with
$\Ext^i_P(-\varpi_2,
-2\varpi_2)\simeq
\Ext^i_P(\varepsilon,
-\varpi_2)=0$ $\forall i\in\bbN$ 
by Bott.
Then
$\Ext^1_P(-\varpi_2,
M^2)
\simeq
\Ext^1_P(-\varpi_2,
\nabla^P(3\varpi_1-3\varpi_2))
\simeq
\Ext^1_P(\varepsilon,
\nabla^P(3\varpi_1-2\varpi_2))
\simeq\Bbbk$
as in (2), and hence the extension $M^3$
is unique up to isomorphism.

We verify finally that the extension
$M^4=\soc^1_{3,e}$
of
$\nabla^P(2\varpi_1-2\varpi_2)$ by $M^3$ is also non-split and uniquely.
Just suppose it split.
As $-2\varpi_1$ is the lowest weight of
$\nabla^P(2\varpi_1-2\varpi_2)$, there would be a $G_1B$-submodule $M''$
of $\soc^3\hat\nabla_P(\varepsilon)$ containing  $\soc^2\hat\nabla_P(\varepsilon)$ to form an exact sequence
$
0\to\soc^2\hat\nabla_P(\varepsilon)\to M''\to
-2p\varpi_1\to 0$,
which would induce an exact sequence of $G$-modules
\begin{multline}
0\to
\ind_{G_1B}^B((2p\varpi_1)\otimes\soc^2\hat\nabla_P(\varepsilon))
\to
\ind_{G_1B}^G((2p\varpi_1)\otimes
M'')\to
\ind_{G_1B}^G(\varepsilon)
\\
\to
\rR^1\ind_{G_1B}^B((2p\varpi_1)\otimes\soc^2\hat\nabla_P(\varepsilon)).
\end{multline}
There are isomorphisms of $G$-modules
$\ind_{G_1B}^G(\varepsilon)\simeq
\nabla(\varepsilon)^{[1]}=L(e)$,
\begin{align*}
\ind_{G_1B}^G
&
((2p\varpi_1)\otimes\soc^2\hat\nabla_P(\varepsilon))
=
\ind_{G_1B}^G((2p\varpi_1)\otimes
\begin{tabular}{|c|}
\hline
$L(s_2)\otimes
(-\varpi_2)^{[1]}$
\\
\hline
$\varepsilon$
\\
\hline
\end{tabular}
)
\\
&\simeq
\ind_{G_1B}^G(
\begin{tabular}{|c|}
\hline
$L(s_2)\otimes
(2\varpi_1-\varpi_2)^{[1]}$
\\
\hline
$(2\varpi_1)^{[1]}$
\\
\hline
\end{tabular}
)
\\
&\simeq
\nabla(2\varpi_1)^{[1]}
\quad\text{by Bott}
\\
&=
L(2\varpi_1)^{[1]}
\quad\text{by the linkage principle
\cite[II.6.17]{J}
under the assumption $p\geq11$}
\\
&\simeq
L(2p\varpi_1).
\end{align*}
Likewise
\begin{align*}
\rR^1\ind_{G_1B}^G((2p\varpi_1)\otimes\soc^2\hat\nabla_P(\varepsilon))
&=
\rR^1\ind_{G_1B}^G((2p\varpi_1)\otimes
\begin{tabular}{|c|}
\hline
$L(s_2)\otimes
(-\varpi_2)^{[1]}$
\\
\hline
$\varepsilon$
\\
\hline
\end{tabular}
)
\\
&\simeq
\rR^1\ind_{G_1B}^G(
\begin{tabular}{|c|}
\hline
$L(s_2)\otimes
(2\varpi_1-\varpi_2)^{[1]}$
\\
\hline
$(2\varpi_1)^{[1]}$
\\
\hline
\end{tabular}
)
=0.
\end{align*}
Thus the exact sequence (3) reads
\[
0\to
L(2p\varpi_1)
\to
\ind_{G_1B}^G((2p\varpi_1)\otimes
M'')
\to
L(e)\to 0.
\]
Also, 
\begin{align*}
\ind_{G_1B}^G((2p\varpi_1)\otimes
M'')
&\leq
\ind_{G_1B}^G((2p\varpi_1)\otimes
\hat\nabla_P(\varepsilon))
\simeq
\ind_{G_1B}^G(\hat\nabla_P(2p\varpi_1))
\\
&\leq
\ind_{G_1B}^G(\hat\nabla(2p\varpi_1))
\simeq
\nabla(2p\varpi_1).
\end{align*}
As $L(2p\varpi_1))
$ is the $G$-socle of $\nabla(2p\varpi_1)$, we must have $\Ext_G^1(L(e),L(2p\varpi_1))\ne0$.
But the distance between the alcoves containing
$0$ and $2p\varpi_1$ is 12 even, contradicting 
\cite[2.10]{A86}/\cite[C.3]{J}. 
Thus, $\soc_{3,e}^1$ is a nonsplit
$P$-extension of
$\nabla^P(2\varpi_1-2\varpi_2)$ by $M^3$.
We show next that the extension is unique up to isomorphism.
One has for each $i\in\bbN$
\begin{align*}
\Ext_P^i
&
(\nabla^P(2\varpi_1-2\varpi_2), -\varpi_2)
\simeq
\Ext_P^i(\varepsilon,
(-\varpi_2)\otimes\nabla^P(2\varpi_1-2\varpi_2)^*)
\\
&\simeq
\Ext_P^i(\varepsilon,
(-\varpi_2)\otimes\nabla^P(2\varpi_1))
\simeq
\Ext_P^i(\varepsilon,
\nabla^P(2\varpi_1-\varpi_2))
\simeq
\Ext_B^i(\varepsilon,
2\varpi_1-\varpi_2)
=0,
\end{align*}
and
\begin{align*}
\Ext_P^i
&
(\nabla^P(2\varpi_1-2\varpi_2), -2\varpi_2)
\simeq
\Ext_P^i(\varepsilon,
\nabla^P(2\varpi_1-2\varpi_2))
\simeq
\Ext_B^i(\varepsilon,
2\varpi_1-2\varpi_2)
\\
&=
\Ext_B^i(\varepsilon,
s_2\bullet(-\varpi_1))
=0.
\end{align*}
Then
\begin{align*}
\Ext_P^1
&
(\nabla^P(2\varpi_1-2\varpi_2), M^3)
\simeq
\Ext_P^1
(\nabla^P(2\varpi_1-2\varpi_2), \nabla^P(3\varpi_1-3\varpi_2))
\\
&\simeq
\Ext_P^1
(\varepsilon,
\nabla^P(2\varpi_1)\otimes\nabla^P(3\varpi_1-3\varpi_2))
\simeq
\Ext_P^1
(\varepsilon,
\begin{tabular}{|c|}
\hline
$\nabla^P(5\varpi_1-3\varpi_2)$
\\
\hline
$\nabla^P(3\varpi_1-2\varpi_2)$
\\
\hline
$\nabla^P(\varpi_1-\varpi_2)$
\\
\hline
\end{tabular})
\\
&\simeq
\Ext_P^1
(\varepsilon,
\begin{tabular}{|c|}
\hline
$\nabla^P(5\varpi_1-3\varpi_2)$
\\
\hline
$\nabla^P(3\varpi_1-2\varpi_2)$
\\
\hline
\end{tabular})
=
\Ext_P^1
(\varepsilon,
\begin{tabular}{|c|}
\hline
$\nabla^P(s_2\bullet(-\varpi_1+\varpi_2))$
\\
\hline
$\nabla^P(s_2\bullet0)$
\\
\hline
\end{tabular})
\\
&\simeq
\Ext_P^1
(\varepsilon,\nabla^P(s_2\bullet0))
\simeq\Ext^1_B(\varepsilon,s_2\bullet0)
\\
&\simeq
\Mod_G(L(e),\rR^1\ind_B^G(s_2\bullet0))
\simeq
\Mod_G(L(e),L(e))\simeq\Bbbk,
\end{align*}
as desired.

Comparing $\varpi_2\otimes\soc^1_{3,e}$
with 
the
adjoint representation of $P$ on its Lie algebra
$\frp$,
we obtain

\begin{prop}
There are isomorphisms of $P$-modules
\[
\soc^1_{3,e}=\begin{tabular}{|c|}
\hline
$\nabla^P(2\varpi_1-2\varpi_2)$
\\
\hline
$-\varpi_2$
\\
\hline
$
\nabla^{P}(3\varpi_1-3\varpi_2)$
\\
\hline
$-2\varpi_2$
\\
\hline
\end{tabular}
\simeq
\begin{tabular}{|c|}
\hline
$(-\varpi_2)\oplus\nabla^P(2\varpi_1-2\varpi_2)$
\\
\hline
$
\nabla^{P}(3\varpi_1-3\varpi_2)$
\\
\hline
$-2\varpi_2$
\\
\hline
\end{tabular}
\simeq
(-\varpi_2)\otimes\frp.
\]

\end{prop}

\setcounter{equation}{0}
\noindent
(1.4)
{\bf Corollary:}
All multiplicity spaces
$\soc^1_{i,w}$, $i\in[1,6]$, $w\in W^P$, 
are indecomposable as $P$-modules.

\setcounter{equation}{0}
\begin{center}
$2^\circ$
{\bf 
Decomposition of $F_*\cO_{G/P}$}

\end{center}

\setcounter{equation}{0}
\noindent
(2.1)
Put $\cP=G/P$.
For each $P$-module
$M$
let
$\cL_\cP(M)$ denote the 
$G$-equivariant sheaf over $\cP$ associated to
$M$.
Sheafifying the socle series of
$\hat\nabla_P(\varepsilon)$
one obtains a filtration of $F_*\cO_\cP$ with subquotients
$\coprod_{w\in W^P}L(w)\otimes\cL_\cP(\soc_{i,w}^1)$.
We will show that the filtration splits, i.e., the $G_1T$-socle series of $\hat\nabla_P(\varepsilon)$
geometrically splits in the terminology of \cite{DG},
to give a decomposition of
$F_*\cO_\cP$ into the direct sum
$F_*\cO_\cP=\coprod_{i=1}^6
\coprod_{w\in W^P}L(w)\otimes\cL_\cP(\soc_{i,w}^1)$
and that
$\cL_\cP(\soc_{3,e}^1)$
causes an obstruction to the self extension of
$F_*\cO_\cP$:
$\Ext_\cP^1(F_*\cO_\cP, F_*\cO_\cP)\ne0$.
Put
$\cM=\cL_\cP(\soc_{3,e}^1)$.

Let
$\ell$ denote the length function on $W$ with respect to the simple reflections.
For each $w\in W^P=\{
e, s_2, s_1s_2,s_2s_1s_2,
s_1s_2s_1s_2,
w^P\}$
put
$\cE(w)=\cL_\cP(\soc_{\ell(w)+1,w}^1)$
the $G$-equivariant sheaf over $\cP$ associated to the $P$-module $\soc_{\ell(w)+1,w}^1$.
We know from
\cite[3.3]{KY}, see Appendix, that
the $\cE(w)$, $w\in W^P$, form a Karoubian complete strongly
exceptional collection of 
coherent modules over $\cP$ such that
$\forall x,y\in W^P$,
$\Mod_\cP(\cE(x),\cE(y))\ne0$ iff
$x\geq y$ in the Chevalley-Bruhat order.
Thus, in order to show that the socle series is geometrially split, it is enough to show that
\begin{align}
\Ext^1_\cP(\cE(w),\cM)=0
\quad\forall w\in W^P\text{ with }
\ell(w)\geq4,
\\
\Ext^1_\cP(\cM,\cE(w))=0
\quad\forall w\in W^P\text{ with }
\ell(w)\leq2.
\end{align}
Moreover, 
the $G_1T$-socle series on
$\hat\nabla_P(\varepsilon)$
coincides with its radical series
\cite{AbK},
so that
$\soc_6\hat\nabla_P(\varepsilon)$ coincides with the head
$\hat\nabla_P(\varepsilon)/\rad
(\hat\nabla_P(\varepsilon))$
of $\hat\nabla_P(\varepsilon)$.
We know from
\cite[\S5]{17} that the inclusion
$\soc\hat\nabla_P(\varepsilon)\hookrightarrow\hat\nabla_P(\varepsilon)$ and the quotient
$\hat\nabla_P(\varepsilon)\twoheadrightarrow\hat\nabla_P(\varepsilon)/\rad
(\hat\nabla_P(\varepsilon))$
are both geometrically split to yield direct summands 
$L(e)\otimes\cE(e)$ and
$L(w^P)\otimes
\cE(w^P)$
of
$F_*\cO_\cP$, and hence we have only to deal with
$w\in\{s_2s_1s_2, s_1s_2s_1s_2\}$ in (1) and 
$w=s_2$ in (2).
We will actually show 
that all higher extension modules
$\Ext^i_\cP(\cE(s_2s_1s_2),\cM)$
$
\Ext^i_\cP(\cE(s_1s_2s_1s_2),\cM)$, and
$
\Ext^i_\cP(\cM,\cE(s_2))$,
$i>0$,
vanish.

\setcounter{equation}{0}
\noindent
(2.2)
Let us 
start the computations.
Put $\cB=G/B$ and let
$\cL(M)$ denote the sheaf over $\cB$ associated to a $B$-module $M$.
Let $i\in\bbN$.
One has
isomorphisms of $G$-modules
\begin{align*}
\Ext^i_\cP(
\cM,\cE(s_2))
&\simeq
\Ext^i_\cP(
\cL_\cP(\begin{tabular}{|c|}
\hline
$(-\varpi_2)\oplus\nabla^P(2\varpi_1-2\varpi_2)$
\\
\hline
$
\nabla^P(3\varpi_1-3\varpi_2)$
\\
\hline
$-2\varpi_2$
\\
\hline
\end{tabular}),\cL_\cP(-\varpi_2))
\\
&\simeq
\rH^i(\cP,
\cL_\cP(\begin{tabular}{|c|}
\hline
$2\varpi_2$
\\
\hline
$
\nabla^P(3\varpi_1-3\varpi_2)^*$
\\
\hline
$\varpi_2\oplus\nabla^P(2\varpi_1-2\varpi_2)^*$
\\
\hline
\end{tabular}
\otimes
(-\varpi_2)))
\\
&\simeq
\rH^i(\cP,
\cL_\cP(\begin{tabular}{|c|}
\hline
$\varpi_2$
\\
\hline
$
\nabla^P(3\varpi_1-\varpi_2)$
\\
\hline
$\varepsilon\oplus\nabla^P(2\varpi_1-\varpi_2)$
\\
\hline
\end{tabular}))
\\
&\simeq
\delta_{i,0}\{
L(\varpi_2)\oplus
L(e)\}
\quad\text{by the linkage principle}.
\end{align*}

\setcounter{equation}{0}
\noindent
(2.3)
Note that
$\nabla(\varpi_2)$ coincides with the adjoint representation of $G$ on its Lie algebra, and hence that
$\nabla(\varpi_2)/\frp\simeq
\begin{tabular}{|c|}
\hline
$\varpi_2$
\\
\hline
$\nabla^P(3\varpi_1-\varpi_2)$
\\
\hline
\end{tabular}$.
We will frequently make use of
identifications
\begin{align}
\soc_{3,e}^1\simeq
(-\varpi_2)\otimes\ker(\frg\tto\frg/\frp)
\simeq
(-\varpi_2)\otimes\ker(\nabla(\varpi_2)\tto
\begin{tabular}{|c|}
\hline
$\varpi_2$
\\
\hline
$\nabla^P(3\varpi_1-\varpi_2)$
\\
\hline
\end{tabular}).
\end{align}
For each $i\in\bbN$ one has isomorhisms of $G$-modules
\begin{align*}
\Ext^i_\cP
&
(
\cE(s_2s_1s_2),\cM)
\\
&\simeq
\Ext^i_\cP(
\cL_\cP(\nabla^P(\varpi_1-2\varpi_2)),
\cL_\cP((-\varpi_2)\otimes\ker(\nabla(\varpi_2)\twoheadrightarrow
\begin{tabular}{|c|}
\hline
$\varpi_2$
\\
\hline
$\nabla^P(3\varpi_1-\varpi_2)$
\\
\hline
\end{tabular})))
\\
&\simeq
\rH^i(\cP,
\cL_\cP(\nabla^P(\varpi_1)\otimes\ker(\nabla(\varpi_2)\twoheadrightarrow
\begin{tabular}{|c|}
\hline
$\varpi_2$
\\
\hline
$\nabla^P(3\varpi_1-\varpi_2)$
\\
\hline
\end{tabular}))),
\end{align*}
giving rise to a long exact sequence
\begin{multline}
\dots\to
\Ext^i_\cP(
\cE(s_2s_1s_2),\cM)
\to
\rH^i(\cP,
\cL_\cP(\nabla^P(\varpi_1)\otimes\nabla(\varpi_2)))
\to
\\
\rH^{i}(\cP,\cL_\cP(
\nabla^P(\varpi_1)\otimes\begin{tabular}{|c|}
\hline
$\varpi_2$
\\
\hline
$\nabla^P(3\varpi_1-\varpi_2)$
\\
\hline
\end{tabular}))
\to\dots
\end{multline}
with isomorphisms of $G$-modules
\begin{align*}
\rH^{i}(\cP,
&
\cL_\cP(
\nabla^P(\varpi_1)\otimes\begin{tabular}{|c|}
\hline
$\varpi_2$
\\
\hline
$\nabla^P(3\varpi_1-\varpi_2)$
\\
\hline
\end{tabular}))
\simeq
\rH^{i}(\cP,
\cL_\cP(
\begin{tabular}{|c|}
\hline
$\nabla^P(\rho)$
\\
\hline
$\nabla^P(\varpi_1)\otimes
\nabla^P(3\varpi_1-\varpi_2)$
\\
\hline
\end{tabular}
))
\\
&\simeq
\rH^{i}(\cP,
\cL_\cP(
\begin{tabular}{|c|}
\hline
$\nabla^P(\rho)$
\\
\hline
$\nabla^P(4\varpi_1-\varpi_2)$
\\
\hline
$\nabla^P(2\varpi_1)$
\\
\hline
\end{tabular}
)
\simeq
\delta_{i,0}\begin{tabular}{|c|}
\hline
$\nabla(\rho)$
\\
\hline
$\nabla(2\varpi_1)$
\\
\hline
\end{tabular}
\\
&\simeq
\delta_{i,0}\{
L(\rho)\oplus
L(2\varpi_1)\}
\quad\text{by the linkage principle}
\end{align*}
and
\begin{align*}
\rH^i(\cP,
\cL_\cP(\nabla^P(\varpi_1)\otimes\nabla(\varpi_2)))
&\simeq
\rH^i(\cP,\cL_\cP(\nabla^P(\varpi_1)))\otimes\nabla(\varpi_2)
\simeq
\delta_{i,0}\nabla(\varpi_1)\otimes\nabla(\varpi_2).
\end{align*}
Thus the sequence (2) reads as an exact sequence
\begin{multline}
0\to
\Mod_\cP(
\cE(s_2s_1s_2),\cM)
\to
\nabla(\varpi_1)\otimes\nabla(\varpi_2)
\to
L(\rho)\oplus
L(2\varpi_1)
\to
\\
\Ext^1_\cP(
\cE(s_2s_1s_2),\cM)
\to
0\to\dots
\end{multline}
and  
$\Ext^j_\cP(
\cE(s_2s_1s_2),\cM)
=0$ $\forall j\geq2$.
On the other hand,
\begin{align*}
\Ext^i_\cP
&
(
\cE(s_2s_1s_2),\cM)
\simeq
\rH^{i}(\cP,
\cL_\cP(
\nabla^P(\rho)\otimes\begin{tabular}{|c|}
\hline
$-\varpi_2\oplus\nabla^P(2\varpi_1-2\varpi_2)$
\\
\hline
$\nabla^P(3\varpi_1-3\varpi_2)$
\\
\hline
$-2\varpi_2$
\\
\hline
\end{tabular}))
\\
&\simeq
\rH^{i}(\cP,
\cL_\cP(
\begin{tabular}{|c|}
\hline
$\nabla^P(\varpi_1)\oplus
\nabla^P(\rho)\otimes\nabla^P(2\varpi_1-2\varpi_2)$
\\
\hline
$\nabla^P(\rho)\otimes\nabla^P(3\varpi_1-3\varpi_2)$
\\
\hline
$\nabla^P(\varpi_1-\varpi_2)$
\\
\hline
\end{tabular}))
\\
&\simeq
\rH^{i}(\cP,
\cL_\cP(
\begin{tabular}{|c|}
\hline
$\nabla^P(\varpi_1)$
\\
\hline
$\nabla^P(3\varpi_1-\varpi_2)$
\\
\hline
$\nabla^P(\varpi_1)$
\\
\hline
$\nabla^P(4\varpi_1-2\varpi_2)$
\\
\hline
$\nabla^P(2\varpi_1-\varpi_2)$
\\
\hline
$\nabla^P(\varpi_1-\varpi_2)$
\\
\hline
\end{tabular}))
\simeq
\rH^{i}(\cP,
\cL_\cP(
\begin{tabular}{|c|}
\hline
$\nabla^P(\varpi_1)^{\oplus_2}$
\\
\hline
$\nabla^P(s_2\bullet\varpi_1)$
\\
\hline
\end{tabular})),
\end{align*}
which induces another exact sequence
\begin{multline}
0\to
\Mod_\cP(
\cE(s_2s_1s_2),\cM)
\to
L(\varpi_1)^{\oplus_2}
\to
L(\varpi_1)
\to
\Ext^1_\cP(
\cE(s_2s_1s_2),\cM)
\to
0.
\end{multline}
Comparing with (3),
we must have 
\[
\Ext^i_\cP(
\cE(s_2s_1s_2),\cM)\simeq
\delta_{i,0}\nabla(\varpi_1).
\]

\setcounter{equation}{0}
\noindent
(2.4)
For each $i\in\bbN$ one has isomorhisms of $G$-modules
\begin{align*}
\Ext^i_\cP
&
(
\cE(s_1s_2s_1s_2),\cM)
\simeq
\Ext^i_\cP(
\cL_\cP((-2\varpi_2)\otimes(\nabla(\varpi_1)/\nabla^P(\varpi_1-\varpi_2))),
\cM)
\\
&\simeq
\Ext^i_\cP(
\cL_\cP(
\begin{tabular}{|c|}
\hline
$\nabla^P(\varpi_1-2\varpi_2)$
\\
\hline
$\nabla^P(2\varpi_1-3\varpi_2)$
\\
\hline
\end{tabular}
),
\cM)
\simeq
\Ext^i_\cP(
\begin{tabular}{|c|}
\hline
$\cE(s_2s_1s_2)$
\\
\hline
$\cL_\cP(\nabla^P(2\varpi_1-3\varpi_2))$
\\
\hline
\end{tabular}
),
\cM),
\end{align*}
giving rise to a long exact sequence
\begin{multline}
\dots\to
\Ext^i_\cP(
\cE(s_2s_1s_2),\cM)
\to
\Ext^i_\cP(\cE(s_1s_2s_1s_2),\cM)
\to
\\
\Ext^i_\cP(\cL_\cP(
\nabla^P(2\varpi_1-3\varpi_2)),\cM)
\to\dots
\end{multline}
with 
$\Ext^i_\cP(
\cE(s_2s_1s_2),\cM)\simeq\delta_{i,0}\nabla(\varpi_1)$
by
(2.3).
Thus the sequence (1)
reads as
\begin{multline}
0\to
\nabla(\varpi_1)
\to
\Mod_\cP(
\cE(s_1s_2s_1s_2),\cM)
\to
\Mod_\cP(\cL_\cP(
\nabla^P(2\varpi_1-3\varpi_2)),\cM)
\to0
\end{multline}
and gives isomorphisms for each
$j\geq1$
\begin{align}
\Ext^j_\cP
&
(
\cE(s_1s_2s_1s_2),\cM)
\simeq
\Ext^j_\cP(\cL_\cP(
\nabla^P(2\varpi_1-3\varpi_2)),\cM)
\\
\notag&\simeq
\rH^j(\cP,
\cL_\cP(
\nabla^P(2\varpi_1)\otimes\ker(\nabla(\varpi_2)\twoheadrightarrow
\begin{tabular}{|c|}
\hline
$\varpi_2$
\\
\hline
$
\nabla^P(3\varpi_1-\varpi_2)$
\\
\hline
\end{tabular}
))).
\end{align}
There is a long exact sequence
\begin{multline}
\dots\to
\Ext^i_\cP(\cL_\cP(
\nabla^P(2\varpi_1-3\varpi_2)),\cM)
\to
\rH^i(\cP,\cL_\cP(\nabla^P(2\varpi_1)\otimes\nabla(\varpi_2)))
\to
\\
\rH^i(\cP,\cL_\cP(\nabla^P(2\varpi_1)\otimes\begin{tabular}{|c|}
\hline
$\varpi_2$
\\
\hline
$
\nabla^P(3\varpi_1-\varpi_2)$
\\
\hline
\end{tabular}
))
\to\dots
\end{multline}
with 
$\rH^i(\cP,\cL_\cP(\nabla^P(2\varpi_1)\otimes\nabla(\varpi_2)))
\simeq
\delta_{i,0}\nabla(2\varpi_1)\otimes\nabla(\varpi_2) 
$
and also
\begin{align*}
\rH^i
&(\cP,\cL_\cP(\nabla^P(2\varpi_1)\otimes\begin{tabular}{|c|}
\hline
$\varpi_2$
\\
\hline
$
\nabla^P(3\varpi_1-\varpi_2)$
\\
\hline
\end{tabular}
))
\simeq
\rH^i
(\cP,\cL_\cP(\begin{tabular}{|c|}
\hline
$\nabla^P(2\varpi_1+\varpi_2)$
\\
\hline
$\nabla^P(2\varpi_1)\otimes
\nabla^P(3\varpi_1-\varpi_2)$
\\
\hline
\end{tabular}
))
\\
&\simeq
\rH^i
(\cP,\cL_\cP(\begin{tabular}{|c|}
\hline
$\nabla^P(2\varpi_1+\varpi_2)$
\\
\hline
$
\nabla^P(5\varpi_1-\varpi_2)$
\\
\hline
$\nabla^P(3\varpi_1)$
\\
\hline
$\nabla^P(\rho)$
\\
\hline
\end{tabular}
))
\simeq
\delta_{i,0}
\begin{tabular}{|c|}
\hline
$\nabla(2\varpi_1+\varpi_2)$
\\
\hline
$\nabla(3\varpi_1)$
\\
\hline
$\nabla(\rho)$
\\
\hline
\end{tabular}.
\end{align*}
Thus the sequence (4) reads as an exact sequence
\begin{multline}
0\to
\Mod_\cP(
\cL_\cP(
\nabla^P(2\varpi_1-3\varpi_2)),\cM)
\to
\nabla(2\varpi_1)\otimes\nabla(\varpi_2)
\to
\begin{tabular}{|c|}
\hline
$\nabla(2\varpi_1+\varpi_2)$
\\
\hline
$\nabla(3\varpi_1)$
\\
\hline
$\nabla(\rho)$
\\
\hline
\end{tabular}\to
\\
\Ext^1_\cP(
\cL_\cP(
\nabla^P(2\varpi_1-3\varpi_2)),\cM)
\to
0,
\end{multline}
and  
gives
$\Ext^j_\cP(
\nabla^P(2\varpi_1-3\varpi_2)),\cM)
=0$ $\forall j\geq2$.

On the other hand,
\begin{align*}
\Ext^i_\cP
&
(
\nabla^P(2\varpi_1-3\varpi_2)),\cM)
\simeq
\rH^{i}(\cP,
\cL_\cP(
\nabla^P(2\varpi_1+\varpi_2)\otimes
\begin{tabular}{|c|}
\hline
$-\varpi_2\oplus\nabla^P(2\varpi_1-2\varpi_2)$
\\
\hline
$\nabla^P(3\varpi_1-3\varpi_2)$
\\
\hline
$-2\varpi_2$
\\
\hline
\end{tabular}))
\\
&\simeq
\rH^{i}(\cP,
\cL_\cP(
\begin{tabular}{|c|}
\hline
$\nabla^P(2\varpi_1)$
\\
\hline
$\nabla^P(4\varpi_1-\varpi_2)$
\\
\hline
$\nabla^P(2\varpi_1)$
\\
\hline
$\nabla^P(\varpi_2)$
\\
\hline
$\nabla^P(5\varpi_1-2\varpi_2)$
\\
\hline
$\nabla^P(3\varpi_1-\varpi_2)$
\\
\hline
$\nabla^P(\varpi_1)$
\\
\hline
\end{tabular}))
\simeq
\rH^{i}(\cP,
\cL_\cP(
\begin{tabular}{|c|}
\hline
$\nabla^P(2\varpi_1)^{\oplus_2}$
\\
\hline
$\nabla^P(\varpi_2)$
\\
\hline
$\nabla^P(s_2\bullet2\varpi_1)$
\\
\hline
$\nabla^P(\varpi_1)$
\\
\hline
\end{tabular}))
\end{align*}
which induces an exact sequence
\begin{multline}
0\to
L(\varpi_1)\to
\Mod_\cP(
\cL_\cP(\nabla^P(2\varpi_1-3\varpi_2)),\cM)
\to
\begin{tabular}{|c|}
\hline
$\nabla(2\varpi_1)^{\oplus_2}$
\\
\hline
$\nabla(\varpi_2)$
\\
\hline
\end{tabular}
\to
L(2\varpi_1)
\to
\\
\Ext^1_\cP(
\cL_\cP(\nabla^P(2\varpi_1-3\varpi_2),\cM)
\to
0.
\end{multline}
Then, together with (5), we must have
\[
\Ext^i_\cP(
\cL_\cP(\nabla^P(2\varpi_1-3\varpi_2),\cM)
\simeq
\delta_{i,0}
\{
L(\varpi_1)\oplus
L(2\varpi_1)\oplus
L(\varpi_2)
\}.
\]
It now follows from (2) and (3) that
\[
\Ext^i_\cP(
\cE(s_1s_2s_1s_2),\cM)\simeq
\delta_{i,0}\{
L(\varpi_1)^{\oplus_2}\oplus
L(2\varpi_1)\oplus
L(\varpi_2)
\}.
\]

\setcounter{equation}{0}
\noindent
(2.5)
As $\nabla(\varpi_1)$ is simple, $\nabla(\varpi_1)$ is also a Weyl module of highest weight $\varpi_1$, and
there is a closed imbedding
$i:\cP\to\bbP(\nabla(\varpi_1))$
such that
$i^*(\cO_{\bbP(\nabla(\varpi_1))}(1))\simeq
\cL_\cP(\varpi_2)$
\cite[II.8.5]{J}.
For an $\cO_\cP$-module
$\cF$ and $n\in\bbZ$
let us abbreviate
$\cF\otimes_\cP\cL_\cP(n\varpi_2)$
as
$\cF(n)$.
Then
$\cM\simeq
\cL_\cP(\frp)\otimes_\cP\cL_\cP(-\varpi_2)=\cL_\cP(\frp)(-1)$.
We have obtained

\begin{thm}
Assume $p\geq11$.
One has a decomposition
\begin{align*}
F_*\cO_\cP
&\simeq
\cE(e)\otimes
L(e)
\oplus
\cE(s_2)\otimes\{L(s_2)\oplus
L(s_1s_2s_1s_2)\}
\oplus
\cE(s_1s_2)\otimes
L(s_1s_2)
\\
&\hspace{1cm}\oplus
\cL_\cP(\frp)(-1)\otimes
L(e)
\oplus
\cE(s_2s_1s_2)\otimes
L(s_2s_1s_2)
\oplus
\cE(s_1s_2s_1s_2)\otimes
L(s_1s_2s_1s_2)
\\
&\hspace{1cm}\oplus
\cE(s_2s_1s_2s_1s_2)\otimes
\{L(s_2s_1s_2s_1s_2)\oplus
L(s_1s_2)\}.
\end{align*}
\end{thm}

\setcounter{equation}{0}
\begin{center}
$3^\circ$
{\bf 
Extensions}

\end{center}

\setcounter{equation}{0}
\noindent
(3.1)
Let $i\in\bbN$.
One has 
\begin{multline*}
\Ext^i_\cP(\cM,\cE(s_2s_1s_2))
\simeq
\\
\Ext^i_\cP(\cL_\cP((-\varpi_2)\otimes\ker(\nabla(\varpi_2)\twoheadrightarrow
\begin{tabular}{|c|}
\hline
$\varpi_2$
\\
\hline
$
\nabla^P(3\varpi_1-\varpi_2)$
\\
\hline
\end{tabular})),\cL_\cP(
\nabla^P(\varpi_1-2\varpi_2))),
\end{multline*}
which gives a long exact sequence
\begin{multline*}
\dots\to
\Ext^i_\cP(\cL_\cP(
\begin{tabular}{|c|}
\hline
$\varepsilon$
\\
\hline
$
\nabla^P(3\varpi_1-2\varpi_2)$
\\
\hline
\end{tabular}),
\cL_\cP(
\nabla^P(\varpi_1-2\varpi_2)))
\to
\\
\Ext^i_\cP(\cL_\cP(
(-\varpi_2)\otimes\nabla(\varpi_2)),
\cL_\cP(
\nabla^P(\varpi_1-2\varpi_2)))
\to
\Ext^i_\cP(\cM,\cE(s_2s_1s_2))
\to\dots
\end{multline*}
with
\begin{align*}
\Ext^i_\cP
&
(\cL_\cP(
(-\varpi_2)\otimes\nabla(\varpi_2)),
\cL_\cP(
\nabla^P(\varpi_1-2\varpi_2)))
\simeq
\nabla(\varpi_2)^*\otimes
\rH^i(\cP,
\cL_\cP(
\nabla^P(\varpi_1-\varpi_2)))
\\
&\simeq
\nabla(\varpi_2)^*\otimes
\rH^i(\cB,
\cL(
\varpi_1-\varpi_2))
=0.
\end{align*}
There follow isomorphisms
\begin{align*}
\Ext^i_\cP
&
(\cM,\cE(s_2s_1s_2))
\simeq
\Ext^{i+1}_\cP(\cL_\cP(
\begin{tabular}{|c|}
\hline
$\varepsilon$
\\
\hline
$
\nabla^P(3\varpi_1-2\varpi_2)$
\\
\hline
\end{tabular}),
\cL_\cP(
\nabla^P(\varpi_1-2\varpi_2)))
\\
&\simeq
\rH^{i+1}(\cP,
\cL_\cP(
\begin{tabular}{|c|}
\hline
$
\nabla^P(3\varpi_1-2\varpi_2)^*$
\\
\hline
$\varepsilon$
\\
\hline
\end{tabular}\otimes
\nabla^P(\varpi_1-2\varpi_2)))
\\
&\simeq
\rH^{i+1}(\cP,
\cL_\cP(
\begin{tabular}{|c|}
\hline
$
\nabla^P(4\varpi_1-3\varpi_2)$
\\
\hline
$
\nabla^P(2\varpi_1-2\varpi_2)$
\\
\hline
$\nabla^P(\varpi_1-2\varpi_2)$
\\
\hline
\end{tabular}
))
=
\rH^{i+1}(\cP,
\cL_\cP(
\begin{tabular}{|c|}
\hline
$
\nabla^P(s_2s_1\bullet0)$
\\
\hline
$
\nabla^P(s_2\bullet(-\varpi_1)$
\\
\hline
$\nabla^P(s_2s_1\bullet(-\varpi_2))$
\\
\hline
\end{tabular}
))
\\
&\simeq
\rH^{i+1}(\cB,
\cL(
s_2s_1\bullet0))
\simeq\delta_{i+1,2}L(e)=
\delta_{i,1}L(e).
\end{align*}

\setcounter{equation}{0}
\noindent
(3.2)
Together with (2.5) we find
\begin{thm}
Assume $p\geq11$.
One has
$\Ext^1_\cP(F_*\cO_\cP,F_*\cO_\cP)\ne0$.

\end{thm}

\setcounter{equation}{0}
\noindent
(3.3)
Let $\bar\cD_\cP^{(0)}=\cMod_{\cO_\cP^{(1)}}(\cO_\cP,\cO_\cP)$ be the sheaf of rings of small differential operators on $\cP$
with
$\cO_\cP^{(1)}$ denoting the sheaf consisting of the $p$-th powers of the elements of $\cO_\cP$.
This is
the first term of the $p$-filtration 
\cite{Haa}
of 
the sheaf 
of rings of differential operators 
\cite{EGA}, 
and is a central reducion of the sheaf 
of rings of arithmetic differential operators of level 0
\cite{B} 
which is called the sheaf of rings of cristaline 
differential operators 
in \cite{BMR}.

\begin{cor}
Assume $p\geq11$.
One has
$\rH^1(\cP,\bar\cD_\cP^{(0)})\ne0$.

\end{cor}

\setcounter{equation}{0}
\noindent(3.4)
With a little more efforts one can also show 
\begin{prop}
Assume $p\geq11$.
For each $i\in\bbN$ one has
\[
\Ext^i_\cP(\cM,\cM)
\simeq
\Ext^i_\cP(\cL_\cP(\frp),\cL_\cP(\frp))
\simeq
\begin{cases}
L(e)
&\text{if $i=0$},
\\
L(\varpi_1)
&\text{if $i=1$},
\\
0
&\text{else}.
\end{cases}
\]
In particular,
$\cM$ 
is indecomposable as an $\cO_\cP$-modules.

\end{prop}

\pf
Let $i\in\bbN$.
We first show
\begin{equation}
\Ext^\bullet_\cP(\cE(e),\cM)=0=\Ext^\bullet_\cP(\cE(s_2),\cM).
\end{equation}
For 
\begin{align*}
\Ext^i_\cP(\cE(e),\cM)
&\simeq
\rH^i(\cP,\cL_\cP(\begin{tabular}{|c|}
\hline
$(-\varpi_2)\oplus\nabla^P(2\varpi_1-2\varpi_2)$
\\
\hline
$
\nabla^P(3\varpi_1-3\varpi_2)$
\\
\hline
$-2\varpi_2$
\\
\hline
\end{tabular}))
\\
&=
\rH^i(\cP,\cL_\cP(\begin{tabular}{|c|}
\hline
$(-\varpi_2)\oplus\nabla^P(s_2\bullet(-\varpi_1))$
\\
\hline
$
\nabla^P(s_2s_1\bullet(\varpi_1-\varpi_2))$
\\
\hline
$s_2s_1s_2s_1\bullet(-\varpi_2)$
\\
\hline
\end{tabular}))
=0.
\end{align*}
Likewise,
\begin{align*}
\Ext^i_\cP(\cE(s_2),\cM)
&=
\Ext^i_\cP(\cL_\cP(-\varpi_2),\cM)\simeq
\rH^i(\cP,\cL_\cP(\begin{tabular}{|c|}
\hline
$\varepsilon\oplus\nabla^P(2\varpi_1-\varpi_2)$
\\
\hline
$
\nabla^P(3\varpi_1-2\varpi_2)$
\\
\hline
$-\varpi_2$
\\
\hline
\end{tabular}))
\\
&=
\rH^i(\cP,\cL_\cP(\begin{tabular}{|c|}
\hline
$\varepsilon$
\\
\hline
$
\nabla^P(s_2\bullet0)$
\\
\hline
\end{tabular})),
\end{align*}
which gives 
$\Ext^j_\cP(\cE(s_2),\cM)=0$
$\forall j\geq2$, and an exact sequence
\[
0\to
\Mod_\cP(\cE(s_2),\cM)
\to
L(e)\to
L(e)
\to
\Ext^1_\cP(\cE(s_2),\cM)
\to0.
\]
On the other hand,
\begin{align*}
\Ext^i_\cP(\cE(s_2),\cM)
&=
\Ext^i_\cP(\cL_\cP(-\varpi_2),
\cL_\cP((-\varpi_2)\otimes\ker(\nabla(\varpi_2)\tto
\begin{tabular}{|c|}
\hline
$\varpi_2$
\\
\hline
$\nabla^P(3\varpi_1-\varpi_2)$
\\
\hline
\end{tabular}))
)
\\
&\simeq
\rH^i(\cP,
\cL_\cP(\ker(\nabla(\varpi_2)\tto
\begin{tabular}{|c|}
\hline
$\varpi_2$
\\
\hline
$\nabla^P(3\varpi_1-\varpi_2)$
\\
\hline
\end{tabular}))
),
\end{align*}
which gives
$\Mod_\cP(\cE(s_2),\cM)
\leq
\nabla(\varpi_2)=L(\varpi_2)$.
We must then have
$
\Mod_\cP(\cE(s_2),\cM)=0=
\Ext^1_\cP(\cE(s_2),\cM)$ also.

Now,
$
\Ext^i_\cP(\cM,\cM)=
\Ext^i_\cP(\cL_\cP((-\varpi_2)\otimes\ker(\nabla(\varpi_2)\tto
\begin{tabular}{|c|}
\hline
$\varpi_2$
\\
\hline
$\nabla^P(3\varpi_1-\varpi_2)$
\\
\hline
\end{tabular}))
,\cM)$ gives rise to a long exact sequence of $G$-modules
\begin{multline*}
\dots\to
\Ext^i_\cP(\cL_\cP(\begin{tabular}{|c|}
\hline
$\varepsilon$
\\
\hline
$\nabla^P(3\varpi_1-2\varpi_2)$
\\
\hline
\end{tabular}))
,\cM)
\to
\Ext^i_\cP(\cL_\cP((-\varpi_2)\otimes\nabla(\varpi_2))
,\cM)
\to
\\
\Ext^i_\cP(\cM,\cM)\to\dots
\end{multline*}
with
\begin{align*}
\Ext^i_\cP(\cL_\cP((-\varpi_2)\otimes\nabla(\varpi_2))
,\cM)
&
\simeq
\nabla(\varpi_2)^*\otimes
\Ext^i_\cP(\cL_\cP(-\varpi_2)
,\cM)
\\
&=
\nabla(\varpi_2)^*\otimes
\Ext^i_\cP(\cE(s_2)
,\cM)
\\
&=0
\quad\text{by (1)},
\end{align*}
and hence
$\Ext^i_\cP(\cM,\cM)\simeq
\Ext^{i+1}_\cP(\cL_\cP(\begin{tabular}{|c|}
\hline
$\varepsilon$
\\
\hline
$\nabla^P(3\varpi_1-2\varpi_2)$
\\
\hline
\end{tabular}))
,\cM)
$.
There arises then a long exact sequence of $G$-modules
\begin{multline*}
\dots\to
\Ext^{i+1}_\cP(\cE(e)
,\cM)
\to
\Ext^i_\cP(\cM,\cM)
\to
\Ext^{i+1}_\cP(\cL_\cP(\nabla^P(3\varpi_1-2\varpi_2)),\cM)\to\dots
\end{multline*}
with
$\Ext^{\bullet}_\cP(\cE(e)
,\cM)=0$ by (1), and hence
\begin{align}
\Ext^i_\cP
&
(\cM,\cM)
\simeq
\Ext^{i+1}_\cP(\cL_\cP(\nabla^P(3\varpi_1-2\varpi_2)),\cM)
\\
\notag&\simeq
\rH^{i+1}(\cP,
\cL_\cP(\nabla^P(3\varpi_1-2\varpi_2))\otimes\ker(\nabla(\varpi_2)\tto
\begin{tabular}{|c|}
\hline
$\varpi_2$
\\
\hline
$\nabla^P(3\varpi_1-\varpi_2)$
\\
\hline
\end{tabular}))).
\end{align}
One obtains then another long exact sequence
\begin{multline*}
\dots\to
\Ext^{i}_\cP(\cM
,\cM)
\to
\rH^{i+1}(\cP,\cL_\cP(\nabla^P(3\varpi_1-2\varpi_2))\otimes\nabla(\varpi_2)))
\to
\\
\rH^{i+1}(\cP,
\cL_\cP(\nabla^P(3\varpi_1-2\varpi_2)\otimes
\begin{tabular}{|c|}
\hline
$\varpi_2$
\\
\hline
$\nabla^P(3\varpi_1-\varpi_2)$
\\
\hline
\end{tabular}))\to\dots
\end{multline*}
with
\begin{align*}
\rH^{i+1}
&
(\cP,\cL_\cP(\nabla^P(3\varpi_1-2\varpi_2))\otimes\nabla(\varpi_2)))
\simeq
\rH^{i+1}
(\cP,\cL_\cP(\nabla^P(3\varpi_1-2\varpi_2)))\otimes\nabla(\varpi_2)
\\
&\simeq
\rH^{i+1}
(\cB,\cL(s_2\bullet0)))\otimes\nabla(\varpi_2)
\simeq
\delta_{i+1,1}L(e)\otimes\nabla(\varpi_2)
\simeq
\delta_{i,0}\nabla(\varpi_2).
\end{align*} 
Thus the exact sequence above 
reads
\begin{multline}
0
\to
\rH^0(\cP,
\cL_\cP(\nabla^P(3\varpi_1-2\varpi_2)\otimes
\nabla^P(3\varpi_1-\varpi_2)))
\to
\Mod_\cP(\cM
,\cM)
\to
\nabla(\varpi_2)
\to
\\
\rH^{1}(\cP,
\cL_\cP(\nabla^P(3\varpi_1-2\varpi_2)\otimes
\nabla^P(3\varpi_1-\varpi_2)))
\to
\Ext^1_\cP(\cM
,\cM)\to0
\end{multline}
and isomorphisms 
$\Ext^j_\cP
(\cM
,\cM)\simeq
\rH^{j}(\cP,
\cL_\cP(\nabla^P(3\varpi_1-2\varpi_2)\otimes
\nabla^P(3\varpi_1-\varpi_2)))$
$\forall j\geq2$.
Now,
\begin{align*}
\rH^{i}
&
(\cP,
\cL_\cP(\nabla^P(3\varpi_1-2\varpi_2)\otimes
\nabla^P(3\varpi_1-\varpi_2)))
\simeq
\rH^{i
}(\cP,
\cL_\cP(
\begin{tabular}{|c|}
\hline
$\nabla^P(6\varpi_1-3\varpi_2)$
\\
\hline
$\nabla^P(4\varpi_1-2\varpi_2)$
\\
\hline
$\nabla^P(2\varpi_1-\varpi_2)$
\\
\hline
$\varepsilon$
\\
\hline
\end{tabular}))
\\
&=
\rH^{i}(\cP,
\cL_\cP(
\begin{tabular}{|c|}
\hline
$\nabla^P(s_2\bullet\varpi_2)$
\\
\hline
$\nabla^P(s_2\bullet\varpi_1)$
\\
\hline
$\varepsilon$
\\
\hline
\end{tabular}))
\simeq
\begin{cases}
L(e)
&\text{if $i=0$},
\\
L(\varpi_2)\oplus
L(\varpi_1)
&\text{if $i=1$
},
\\
0
&\text{else},
\end{cases}
\end{align*}
and hence the sequence (3) reads
\[
0
\to
L(e)
\to
\Mod_\cP(\cM
,\cM)
\to
L(\varpi_2)
\to
L(\varpi_2)\oplus
L(\varpi_1)
\to
\Ext^1_\cP(\cM
,\cM)\to0.
\]

On the other hand, 
$\Ext^i_\cP
(\cM,\cM)
\simeq
\rH^{i+1}(\cP,
\cL_\cP(\nabla^P(3\varpi_1-2\varpi_2)\otimes\frp))
$
from (2)
with
$\nabla^P(3\varpi_1-2\varpi_2)\otimes\frp\simeq
\nabla^P(3\varpi_1-2\varpi_2)\otimes
\begin{tabular}{|c|}
\hline
$\nabla^P(2\varpi_1-\varpi_2)$
\\
\hline
$\varepsilon$
\\
\hline
$\nabla^P(3\varpi_1-2\varpi_2)$
\\
\hline
$-\varpi_2$
\\
\hline
\end{tabular}
$
having a $P$-module
filtration
whose subquotients
are
$\nabla^P(5\varpi_1-3\varpi_2)=
\nabla^P(s_2\bullet(-\varpi_1+\varpi_2))$,
$\nabla^P(3\varpi_1-2\varpi_2)=
\nabla^P(s_2\bullet0)$ twice,
$\nabla^P(\varpi_1-\varpi_2)
$,
$\nabla^P(6\varpi_1-4\varpi_2)
=
\nabla^P(s_2s_1\bullet\varpi_1)$,
$\nabla^P(4\varpi_1-3\varpi_2)
=
\nabla^P(s_2s_1\bullet0)$,
$\nabla^P(2\varpi_1-2\varpi_2)
=
\nabla^P(s_2\bullet(-\varpi_1))$,
$\nabla^P(-\varpi_2)=-\varpi_2$,
$\nabla^P(3\varpi_1-3\varpi_2)
=
\nabla^P(s_2s_1\bullet(\varpi_1-\varpi_2))$.
It follows that
the possible $G$-composition factors of 
$\Ext^i_\cP
(\cM,\cM)
$ are just
$L(e)$ and $L(\varpi_1)$,
and hence the assertion.
\\

\setcounter{equation}{0}
\begin{center}
{\bf 
Appendixes}

\end{center}

\noindent
{\bf A.}
Keep the notation from the main text.
We assume, in particular, that
$p\geq11$.
We will recover from
\cite{KY}
the proof of the fact that
the
$\cE(w)=\cL_\cP(\soc_{\ell(w)+1}^1)$, $w\in W^P$, form a Karoubian complete strongly exceptional sequence 
in the bounded derived category 
$\rD^b(\coh\cP)$
of coherent sheaves on $\cP$
such that
$\Mod_\cP(\cE(x),\cE(y))\ne0$ iff
$x\geq y$ in the Chevalley-Bruat order.
The present parametrization of the sheaves is twisted from the one in
\cite{KY}
by the 
involution
$w_0?w_P$ on $W^P$,
which reverses
the
Chevalley-Bruhat order.
Incorporating 
progress since,
we also
employ \cite{AbK} to replace some ad hoc arguments
using 
\cite{AK89}.

\setcounter{equation}{0}
\noindent
(A.1)
In order to determine the $G_1P$-module structure on
$\hat\nabla_P(\varepsilon)$,
we first
compute its $G_1T$-structure using the formula
\cite[5.2]{AbK}.
As each
$\Mod_{G_1}(L(w),\soc_i \hat\nabla_P(\varepsilon))$
is equipped with a structure of $P$-module,
one readily finds
in the notation of (1.3)
\begin{align}
\soc_1\hat\nabla_P(\varepsilon)&\simeq
L(e)\otimes\varepsilon,
\\
\notag
\soc_2\hat\nabla_P(\varepsilon)&\simeq
L(s_2)\otimes
(-\varpi_2)^{[1]},
\\
\notag
\soc_3\hat\nabla_P(\varepsilon)&\simeq
L(s_1s_2)\otimes
(\soc_{3,s_1s_2}^1)^{[1]}
\oplus
L(s_1s_2s_1s_2)\otimes
(-\varpi_2)^{[1]}
\oplus
L(e)\otimes(\soc_{3,e}^1)^{[1]},
\\
\notag
\soc_4\hat\nabla_P(\varepsilon)&\simeq
L(s_2s_1s_2)\otimes
\nabla^P(\varpi_1-2\varpi_2)^{[1]},
\\
\notag
\soc_5\hat\nabla_P(\varepsilon)&\simeq
L(s_1s_2s_1s_2)\otimes
(\soc_{5,s_1s_2s_1s_2}^1)^{[1]}\oplus
L(s_1s_2)\otimes(-2\varpi_2)^{[1]},
\\
\notag
\soc_6\hat\nabla_P(\varepsilon)&\simeq
L(
w^P)\otimes
(-2\varpi_2)^{[1]},
\end{align}
with
$\soc_{3,s_1s_2}^1=
\begin{tabular}{|c|}
\hline
$\nabla^P(2\varpi_1-2\varpi_2)$
\\
\hline
$\nabla^P(\varpi_1-2\varpi_2)$
\\
\hline
\end{tabular}$
and
$\soc_{5,s_1s_2s_1s_2}^1=
\begin{tabular}{|c|}
\hline
$\nabla^P(\varpi_1-2\varpi_2)$
\\
\hline
$\nabla^P(2\varpi_1-3\varpi_2)$
\\
\hline
\end{tabular}$. 

We show that both of the last two extensions as $P$-modules
are nonsplit and uniquely
to yield 
isomorphisms of $P$-modules
\begin{align}
\soc^1_{3,s_1s_2}
&\simeq
(-\varpi_2)\otimes\ker(\nabla(\varpi_1)\twoheadrightarrow\nabla^P(\varpi_1)),
\\
\notag
\soc^1_{5,s_1s_2s_1s_2}
&\simeq
(-2\varpi_2)\otimes(\nabla(\varpi_1)/\nabla^P(\varpi_1-\varpi_2)).
\end{align}

Just suppose the extension in 
$\soc_{3,s_1s_2}^1$ is split.
Then there would be a 
$G_1B$-submodule $M_1$
of
$\hat\nabla_P(\varepsilon)$
containing 
$\soc^2\hat\nabla_P(\varepsilon)$
such that
$M_1/\soc^2\hat\nabla_P(\varepsilon)\simeq
L(s_1s_2)\otimes(-2\varpi_1)^{[1]}$
as $G_1B$-modules.
It would then induce an exact sequence of
$G$-modules
\[
\ind_{G_1B}^G(M_1\otimes
2p\varpi_1)\to
L(s_1s_2)\to
\rR^1\ind_{G_1B}^G(\soc^2\hat\nabla_P(\varepsilon)\otimes
2p\varpi_1).
\]
But
$\rR^1\ind_{G_1B}^G(\soc^2\hat\nabla_P(\varepsilon)\otimes
2p\varpi_1)$ has no
$G$-composition factor whose $G_1$-part is
$L(s_1s_2)$
while
$\ind_{G_1B}^G(M_1\otimes
2p\varpi_1)
\leq
\ind_{G_1B}^G(\hat\nabla_P(\varepsilon)\otimes
2p\varpi_1)
\leq
\ind_{G_1B}^G(\hat\nabla(\varepsilon)\otimes
2p\varpi_1)
\simeq
\nabla(
2p\varpi_1)$
with
$\nabla(
2p\varpi_1)$ having no composition factor
$L(s_1s_2)$, absurd.
Also,
\begin{align*}
\Ext^1_P
&(\nabla^P(2\varpi_1-2\varpi_2),\nabla^P(\varpi_1-2\varpi_2))
\simeq
\Ext^1_P
(\nabla^P(2\varpi_1),\nabla^P(\varpi_1))
\\
&\simeq
\Ext^1_B
(\nabla^P(2\varpi_1),\varpi_1)
\quad\text{by the Frobenius reciprocity}
\\
&\simeq
\rH^1(B,\varpi_1\otimes
\nabla^P(2\varpi_1)^*)\simeq
\rH^1(B,\varpi_1\otimes
\nabla^P(2\varpi_1-2\varpi_2))
\simeq
\rH^1(B,\begin{tabular}{|c|}
\hline
$3\varpi_1-2\varpi_2$
\\
\hline
$\varpi_1-\varpi_2$
\\
\hline
$-\varpi_1$
\\
\hline
\end{tabular})
\\
&\simeq
\rH^1(B,3\varpi_1-2\varpi_2)
\quad\text{
as
$\rH^\bullet(\cB,\cL(-\varpi_1))=0=
\rH^\bullet(\cB,\cL(\varpi_1-\varpi_2))$
\cite[II.5.4]{J}
}
\\
&=
\rH^1(B,-\alpha_2)
\simeq
\Bbbk.
\end{align*}

Likewise,
just suppose the extension in 
$\soc_{5,s_1s_2s_1s_2}^1$ is split.
Dualizing
$\hat\nabla_P(\varepsilon)$,
by the rigidity of 
$\hat\nabla_P(\varepsilon)^*\simeq
\hat\nabla_P(3(p-1)\varpi_2)$
\cite{AbK}
there would be a $G_1B$-submodule
$M_2$ of
$\soc^2\hat\nabla_P(3(p-1)\varpi_2)$
containing
$\soc\hat\nabla_P(3(p-1)\varpi_2)\simeq
L(w^P)\otimes(2\varpi_2)^{[1]}$
such that
$M_2/\soc\hat\nabla_P(3(p-1)\varpi_2)\simeq
L(s_1s_2s_1s_2)\otimes
(-2\varpi_1+3\varpi_2)^{[1]}$
as $G_1B$-modules.
There would then be an exct sequence of $G$-modules
\[
\ind_{G_1B}^G(M_2\otimes
p(2\varpi_1-3\varpi_2))\to
L(s_1s_2s_1s_2)\to
\rR^1\ind_{G_1B}^G(L(w^P)\otimes
p(2\varpi_1-\varpi_2)).
\]
with
$\rR^1\ind_{G_1B}^G(L(w^P)\otimes
p(2\varpi_1-\varpi_2))\simeq
L(w^P)\otimes
\rR^1\ind_{B}^G(2\varpi_1-\varpi_2)^{[1]}=0$.
But
$\ind_{G_1B}^G(M_2\otimes
p(2\varpi_1-3\varpi_2))
\leq
\ind_{G_1B}^G(\hat\nabla_P(3(p-1)\varpi_2)\otimes
p(2\varpi_1-3\varpi_2))
\leq
\ind_{G_1B}^G(\hat\nabla(3(p-1)\varpi_2)\otimes
p(2\varpi_1-3\varpi_2))
\simeq
\nabla(
2p\varpi_1-3\varpi_2)$
with
$\nabla(
2p\varpi_1-3\varpi_2)$ having no composition factor
$L(s_1s_2s_1s_2)$, absurd.
Also,
\begin{align*}
\Ext^1_P
&(\nabla^P(\varpi_1-2\varpi_2),\nabla^P(2\varpi_1-3\varpi_2))
\simeq
\rH^1(B,(2\varpi_1-\varpi_2)\otimes
\nabla^P(\varpi_1-\varpi_2))
\\
&\simeq
\rH^1(B,\begin{tabular}{|c|}
\hline
$3\varpi_1-2\varpi_2$
\\
\hline
$\varpi_1-\varpi_2$
\\
\hline
\end{tabular})
\simeq
\rH^1(B,3\varpi_1-2\varpi_2)
=
\rH^1(B,-\alpha_2)
\simeq
\Bbbk.
\end{align*}

\begin{prop}
Each $\soc_{\ell(w)+1}^1$, $w\in W^P$, is an indecomposable $P$-module, of 
highest weight $w^{-1}\bullet(w\bullet0)^1$
except for $w=s_1s_2s_1s_2$.
In the last case
$\soc^1_{5,s_1s_2s_1s_2}$ is 
generated by a vector of weight
$\varpi_1-2\varpi_2$.

\end{prop}

\setcounter{equation}{0}
\noindent
(A.2)
We set 
$\cE(w)=\cL_\cP(\soc_{\ell(w)+1}^1)$ for each $w\in W^P$, and determine their mutual extensions
$\Ext^i_\cP(\cE(x),\cE(y))$,
$x,y\in W^P$,
$i\in\bbN$, as $G$-modules.

We start with the computations involving $\cE(e)$. We have
\begin{align}
\Ext^i_\cP(\cE(e),\cE(e))
&\simeq
\rH^i(\cP,\cO_\cP)\simeq
\delta_{i,0}\Bbbk,
\\
\Ext^i_\cP(\cE(e),\cE(s_2))&\simeq
\rH^i(\cP,\cL_\cP(-\varpi_2))\simeq
0,
\\
\Ext^i_\cP(\cE(e),\cE(s_1s_2))&\simeq
\rH^i(\cP,\cL_\cP(\begin{tabular}{|c|}
\hline
$\nabla^P(2\varpi_1-2\varpi_2)$
\\
\hline
$\nabla^P(\varpi_1-2\varpi_2)$
\\
\hline
\end{tabular}))
\simeq
\rH^i(\cB,\cL(\begin{tabular}{|c|}
\hline
$2\varpi_1-2\varpi_2$
\\
\hline
$\varpi_1-2\varpi_2$
\\
\hline
\end{tabular}))
\\
\notag&=
\rH^i(\cB,\cL(\begin{tabular}{|c|}
\hline
$s_2\bullet(-\varpi_1)$
\\
\hline
$s_2s_1\bullet(-\varpi_2)$
\\
\hline
\end{tabular}))
=0
\quad\text{\cite[II.5.5]{J}},
\\
\Ext^i_\cP(\cE(e),\cE(s_2s_1s_2))&\simeq
\rH^i(\cP,\cL_\cP(\nabla^P(\varpi_1-2\varpi_2)))
=
0
\quad\text{as in (3)},
\\
\Ext^i_\cP(\cE(e),\cE(s_1s_2s_1s_2))
&\simeq
\rH^i(\cP,\cL_\cP(\begin{tabular}{|c|}
\hline
$\nabla^P(\varpi_1-2\varpi_2)$
\\
\hline
$\nabla^P(2\varpi_1-3\varpi_2)$
\\
\hline
\end{tabular}))
\simeq
\rH^i(\cB,\cL(\begin{tabular}{|c|}
\hline
$\varpi_1-2\varpi_2$
\\
\hline
$2\varpi_1-3\varpi_2$
\\
\hline
\end{tabular}))
\\
\notag&=
\rH^i(\cB,\cL(\begin{tabular}{|c|}
\hline
$s_2s_1\bullet(-\varpi_2)$
\\
\hline
$s_1s_2s_1\bullet(-\varpi_1)$
\\
\hline
\end{tabular}))
=0,
\\
\Ext^i_\cP(\cE(e),\cE(w^P))
&\simeq
\rH^i(\cP,\cL_\cP(-2\varpi_2))=
\rH^i(\cP,\cL_\cP(s_2s_1s_2s_1\bullet(-\varpi_2))=
0.
\end{align}

\setcounter{equation}{0}
\noindent
(A.3)
We compute next the extensions with
$\cE(s_2)$.
Let
$i\in\bbN$.
One has
\begin{align}
\Ext^i_\cP
&
(\cE(s_2),\cE(e))
\simeq
\rH^i(\cP,\cL_\cP(\varpi_2))\simeq
\delta_{i,0}\nabla(\varpi_2),
\\
\Ext^i_\cP
&
(\cE(s_2),\cE(s_2))
\simeq
\rH^i(\cP,\cO_\cP)\simeq
\delta_{i,0}\Bbbk,
\\
\Ext^i_\cP
&
(\cE(s_2),\cE(s_1s_2))
\simeq
\rH^i(\cP,\cL_\cP(\varpi_2\otimes\begin{tabular}{|c|}
\hline
$\nabla^P(2\varpi_1-2\varpi_2)$
\\
\hline
$\nabla^P(\varpi_1-2\varpi_2)$
\\
\hline
\end{tabular}))
\\
\notag
&\simeq
\rH^i(\cB,\cL(\begin{tabular}{|c|}
\hline
$2\varpi_1-\varpi_2$
\\
\hline
$\varpi_1-\varpi_2$
\\
\hline
\end{tabular}))
=0,
\\
\Ext^i_\cP
&
(\cE(s_2),\cE(s_2s_1s_2))
\simeq
\rH^i(\cP,\cL_\cP(\varpi_2\otimes\nabla^P(\varpi_1-2\varpi_2)))
\\
\notag&\simeq
\rH^i(\cP,\cL_\cP(\nabla^P(\varpi_1-\varpi_2)))
=
0,
\\
\Ext^i_\cP
&
(\cE(s_2),\cE(s_1s_2s_1s_2))
\simeq
\rH^i(\cP,\cL_\cP(\varpi_2\otimes
\begin{tabular}{|c|}
\hline
$\nabla^P(\varpi_1-2\varpi_2)$
\\
\hline
$\nabla^P(2\varpi_1-3\varpi_2)$
\\
\hline
\end{tabular}))
\\
\notag&\simeq
\rH^i(\cB,\cL(\begin{tabular}{|c|}
\hline
$\varpi_1-\varpi_2$
\\
\hline
$2\varpi_1-2\varpi_2$
\\
\hline
\end{tabular}))
\simeq
\rH^i(\cB,\cL(s_2\bullet(-\varpi_1)))
=0,
\\
\Ext^i_\cP&
(\cE(s_2),\cE(w^P))
\simeq
\rH^i(\cP,\cL_\cP(-\varpi_2))=
0.
\end{align}

\setcounter{equation}{0}
\noindent
(A.4)
Let
$i\in\bbN$.
We have
\begin{align}
\Ext^i_\cP
&
(\cE(s_1s_2),\cE(e))
\simeq
\rH^i(\cP,\cL_\cP(\begin{tabular}{|c|}
\hline
$\nabla^P(\varpi_1-2\varpi_2)^*$
\\
\hline
$
\nabla^P(2\varpi_1-2\varpi_2)^*$
\\
\hline
\end{tabular}))
\simeq
\rH^i(\cP,\cL_\cP(\begin{tabular}{|c|}
\hline
$\nabla^P(\rho)$
\\
\hline
$
\nabla^P(2\varpi_1)$
\\
\hline
\end{tabular}))
\\
\notag&\simeq
\delta_{i,0}\{
\nabla(\rho)\oplus\nabla(2\varpi_1)
\}
\quad\text{by the linkage principle
\cite[II.6.17]{J}},
\\
\Ext^i_\cP
&
(\cE(s_1s_2),\cE(s_2))
\simeq
\rH^i(\cP,\cL_\cP((-\varpi_2)\otimes
\begin{tabular}{|c|}
\hline
$\nabla^P(\rho)$
\\
\hline
$
\nabla^P(2\varpi_1)$
\\
\hline
\end{tabular}))
\quad\text{by (1)}
\\
\notag&\simeq
\rH^i(\cP,\cL_\cP(
\begin{tabular}{|c|}
\hline
$\nabla^P(\varpi_1)$
\\
\hline
$
\nabla^P(2\varpi_1-\varpi_2)$
\\
\hline
\end{tabular}))
\simeq\delta_{i,0}\nabla(\varpi_1).
\end{align}

As
$\cE(s_1s_2)\simeq\cL_\cP((-\varpi_2)\otimes\ker(\nabla(\varpi_1)\twoheadrightarrow\nabla^P(\varpi_1)))$
by (A.1.2), there is a long exact sequence
\begin{multline}
0\to
\Mod_\cP(\cL_\cP((-\varpi_2)\otimes\nabla^P(\varpi_1)),\cE(s_1s_2))
\to
\\
\Mod_\cP(\cL_\cP((-\varpi_2)\otimes\nabla(\varpi_1)),\cE(s_1s_2))
\\
\to
\Mod_\cP(\cE(s_1s_2),\cE(s_1s_2))
\to
\dots
\end{multline}
with
\begin{align*}
\Ext^\bullet_\cP
&
(\cL_\cP((-\varpi_2)\otimes\nabla(\varpi_1)),\cE(s_1s_2))
\simeq
\Ext^\bullet_\cP(\cL_\cP(-\varpi_2),\cE(s_1s_2))\otimes
\nabla(\varpi_1)^*
\\
&\hspace{1cm}\text{by the tensor identity as $\nabla(\varpi_1)$ is equipped with a structure of $G$-module}
\\
&
\simeq
\Ext^\bullet_\cP(\cE(s_2),\cE(s_1s_2))\otimes
\nabla(\varpi_1)^*
=0
\quad\text{by
(A.3.3)
}.
\end{align*}
Thus,
$
\Ext^i_\cP(\cE(s_1s_2),\cE(s_1s_2))
\simeq
\Ext^{i+1}_\cP(\cL_\cP((-\varpi_2)\otimes\nabla^P(\varpi_1)),\cE(s_1s_2))$,
the right hand side of which fits into another long exact sequence
\begin{multline*}
0\to
\Mod_\cP(\cL_\cP((-\varpi_2)\otimes\nabla^P(\varpi_1)),\cE(s_1s_2))
\to
\\
\Mod_\cP(\cL_\cP((-\varpi_2)\otimes\nabla^P(\varpi_1)),\cL_\cP((-\varpi_2)\otimes\nabla(\varpi_1)))
\\
\to
\Mod_\cP(\cL_\cP((-\varpi_2)\otimes\nabla^P(\varpi_1)),\cL_\cP((-\varpi_2)\otimes\nabla^P(\varpi_1)))
\to
\dots
\end{multline*}
with
$\Ext^\bullet_\cP(\cL_\cP((-\varpi_2)\otimes\nabla^P(\varpi_1)),\cL_\cP((-\varpi_2)\otimes\nabla(\varpi_1)))
\simeq
\rH^\bullet(\cP,\cL_\cP(\nabla^P(\varpi_1)^*))\otimes
\nabla(\varpi_1)
\simeq
\rH^\bullet(\cP,\cL_\cP(\nabla^P(\varpi_1-\varpi_2)))\otimes
\nabla(\varpi_1)
=0$.
It follows that
\begin{align}
\Ext^i_\cP
&
(\cE(s_1s_2),\cE(s_1s_2))
\simeq
\Ext^i_\cP(\cL_\cP((-\varpi_2)\otimes\nabla^P(\varpi_1)),\cL_\cP((-\varpi_2)\otimes\nabla^P(\varpi_1)))
\\
\notag&\simeq
\rH^i(\cP,
\cL_\cP(\nabla^P(\varpi_1)^*\otimes\nabla^P(\varpi_1)))
\simeq
\rH^i(\cP,
\cL_\cP(\nabla^P(\varpi_1-\varpi_2)\otimes\nabla^P(\varpi_1)))
\\
\notag&\simeq
\rH^i(\cP,\cL_\cP(\begin{tabular}{|c|}
\hline
$\nabla^P(2\varpi_1-\varpi_2)$
\\
\hline
$\nabla^P(\varepsilon)$
\\
\hline
\end{tabular}))
\simeq
\delta_{i,0}\Bbbk.
\end{align}

To find
$\Ext^i_\cP
(\cE(s_1s_2),\cE(s_2s_1s_2))$,
consider the long exact sequence (3) with
$\cE(s_1s_2)$ in the covariant entries
replaced by
$\cE(s_2s_1s_2)$.
As
$\Ext^\bullet_\cP(\cL_\cP((-\varpi_2)\otimes\nabla(\varpi_1)),\cE(s_2s_1s_2))
\simeq
\Ext^\bullet_\cP(\cE(s_2),\cE(s_2s_1s_2))\otimes\nabla(\varpi_1)^*
=0$ by
(A.3.4),
one obtains
\begin{align}
\Ext^i_\cP
&
(\cE(s_1s_2),\cE(s_2s_1s_2))
\simeq
\Ext^{i+1}_\cP(\cL_\cP((-\varpi_2)\otimes\nabla^P(\varpi_1)),\cE(s_2s_1s_2))
\\
\notag&=
\Ext^{i+1}_\cP(\cL_\cP((-\varpi_2)\otimes\nabla^P(\varpi_1)),\cL_\cP(\nabla^P(\varpi_1-2\varpi_2))
\\
\notag&\simeq
\rH^{i+1}(\cP,
\cL_\cP(\nabla^P(\varpi_1-\varpi_2)\otimes\nabla^P(\varpi_1-\varpi_2)))
\\
\notag&\simeq
\rH^{i+1}(\cP,
\cL_\cP(\nabla^P(
\begin{tabular}{|c|}
\hline
$\nabla^P(2\varpi_1-2\varpi_2)$
\\
\hline
$\nabla^P(-\varpi_2)$
\\
\hline
\end{tabular})))
=0
\quad\text{as in (A.2.3)}.
\end{align}
Then
\begin{align*}
\Ext^i_\cP
&
(\cE(s_1s_2),\cE(s_1s_2s_1s_2))
\simeq
\Ext^{i}_\cP(\cE(s_1s_2),\cL_\cP(\begin{tabular}{|c|}
\hline
$\nabla^P(\varpi_1-2\varpi_2)$
\\
\hline
$\nabla^P(2\varpi_1-3\varpi_2)$
\\
\hline
\end{tabular}))
\\
\notag&\simeq
\Ext^{i}_\cP(\cE(s_1s_2),\cL_\cP(\nabla^P(2\varpi_1-3\varpi_2))) \quad\text{by (5)},
\end{align*}
the right hand side of which fits into the long exact sequence
(3)
with
$\cE(s_1s_2)$ in the covariant entries
replaced by
$\cL_\cP(\nabla^P(2\varpi_1-3\varpi_2))$.
As
$\Ext^\bullet_\cP(\cL_\cP((-\varpi_2)\otimes\nabla(\varpi_1)),\cL_\cP(\nabla^P(2\varpi_1-3\varpi_2)))
\simeq
\nabla(\varpi_1)^*\otimes
\rH^\bullet(\cP,
\cL_\cP(\nabla^P(2\varpi_1-2\varpi_2)))=0$, 
one obtains
\begin{align}
\Ext^i_\cP
&
(\cE(s_1s_2),\cE(s_1s_2s_1s_2))
\simeq
\Ext^{i+1}_\cP(\cL_\cP((-\varpi_2)\otimes\nabla^P(\varpi_1)),\cL_\cP(\nabla^P(2\varpi_1-3\varpi_2)))
\\
\notag&\simeq
\rH^{i+1}(\cP,
\cL_\cP(\nabla^P(\varpi_1-\varpi_2)\otimes
\nabla^P(2\varpi_1-2\varpi_2))
\\
\notag&\simeq
\rH^{i+1}(\cP,
\cL_\cP(\nabla^P(
\begin{tabular}{|c|}
\hline
$\nabla^P(3\varpi_1-3\varpi_2)$
\\
\hline
$\nabla^P(\varpi_1-2\varpi_2)$
\\
\hline
\end{tabular})))
\\
\notag&\simeq
\rH^{i+1}(\cP,
\cL_\cP(\nabla^P(
3\varpi_1-3\varpi_2)))
\quad\text{by (A.2.4)}
\\
\notag&=
\rH^{i+1}(\cP,
\cL_\cP(\nabla^P(
s_1s_2\bullet(\varpi_1-\varpi_2)))
=0.
\end{align}

Finally, one has
\begin{align}
\Ext^i_\cP
&
(\cE(s_1s_2),\cE(w^P))
\simeq
\rH^{i}(\cP,
\cL_\cP((-2\varpi_2)\otimes
\begin{tabular}{|c|}
\hline
$\nabla^P(\rho)$
\\
\hline
$
\nabla^P(2\varpi_1)$
\\
\hline
\end{tabular}
))
\quad\text{as in (1)}
\\
\notag&\simeq
\rH^{i}(\cP,
\cL_\cP(
\begin{tabular}{|c|}
\hline
$\nabla^P(\varpi_1-\varpi_2)$
\\
\hline
$\nabla^P(2\varpi_1-2\varpi_2)$
\\
\hline
\end{tabular}))
\simeq
\rH^{i}(\cP,
\cL_\cP(
\nabla^P(2\varpi_1-2\varpi_2)))
\\
\notag&=0
\quad\text{as in (A.2.3)}.
\end{align}

\setcounter{equation}{0}
\noindent
(A.5)
Let
$i\in\bbN$.
We have
\begin{align}
\Ext^i_\cP
&
(\cE(s_2s_1s_2),\cE(e))
\simeq
\rH^i(\cP,\cL_\cP(\nabla^P(\varpi_1-2\varpi_2)^*))
\simeq
\rH^i(\cP,\cL_\cP(\nabla^P(\rho)))
\\
\notag&\simeq
\delta_{i,0}\nabla(\rho),
\\
\Ext^i_\cP
&
(\cE(s_2s_1s_2),\cE(s_2))
\simeq
\rH^i(\cP,\cL_\cP((-\varpi_2)\otimes\nabla^P(\rho)))
\simeq
\delta_{i,0}\nabla(\varpi_1),
\\
\Ext^i_\cP
&
(\cE(s_2s_1s_2),\cE(s_1s_2))
\simeq
\rH^i(\cP,\cL_\cP(\nabla^P(\rho)
\otimes
\begin{tabular}{|c|}
\hline
$\nabla^P(2\varpi_1-2\varpi_2)$
\\
\hline
$\nabla^P(\varpi_1-2\varpi_2)$
\\
\hline
\end{tabular}))
\\
\notag&\simeq
\rH^i(\cP,\cL_\cP(
\begin{tabular}{|c|}
\hline
$\nabla^P(3\varpi_1-\varpi_2)$
\\
\hline
$\nabla^P(\varpi_1)$
\\
\hline
$\nabla^P(2\varpi_1-\varpi_2)$
\\
\hline
$\varepsilon$
\\
\hline
\end{tabular}))
\simeq
\delta_{i,0}
\{
\nabla(\varpi_1)\oplus
\varepsilon
\}
\quad\text{by the linkage principle},
\\
\Ext^i_\cP
&
(\cE(s_2s_1s_2),\cE(s_2s_1s_2))
\simeq
\rH^i(\cP,\cL_\cP(\nabla^P(\rho)
\otimes
\nabla^P(\varpi_1-2\varpi_2)
))
\\
\notag&\simeq
\rH^i(\cP,\cL_\cP(\begin{tabular}{|c|}
\hline
$\nabla^P(2\varpi_1-\varpi_2)$
\\
\hline
$\varepsilon$
\\
\hline
\end{tabular}))
\simeq
\delta_{i,0}\Bbbk.
\end{align}

By (A.1.2)
there is a long exact sequence
\begin{multline}
0\to
\Mod_\cP(\cE(s_2s_1s_2),
\cL_\cP(\nabla^P(\varpi_1-3\varpi_2)))
\\
\to
\Mod_\cP(\cE(s_2s_1s_2),\cL_\cP((-2\varpi_2)\otimes
\nabla(\varpi_1)))
\\
\to
\Mod_\cP(\cE(s_2s_1s_2),\cE(s_1s_2s_1s_2))\to
\dots
\end{multline}
with
\begin{align*}
\Ext^\bullet_\cP
&
(\cE(s_2s_1s_2),\cL_\cP((-2\varpi_2)\otimes
\nabla(\varpi_1)))
\simeq
\rH^i(\cP,\cL_\cP(\nabla^P(\rho)\otimes
(-2\varpi_2)))\otimes
\nabla(\varpi_1)
\\
\notag&\simeq
\rH^i(\cP,\cL_\cP(\nabla^P(
\varpi_1-\varpi_2)))\otimes
\nabla(\varpi_1)
=0,
\end{align*}
and hence
\begin{align}
\Ext^i_\cP
&
(\cE(s_2s_1s_2),\cE(s_1s_2s_1s_2))
\simeq
\Ext^{i+1}_\cP(\cE(s_2s_1s_2),
\cL_\cP(\nabla^P(\varpi_1-3\varpi_2)))
\\
\notag&\simeq
\rH^{i+1}(\cP,
\cP(\nabla^P(\rho)\otimes
\nabla^P(\varpi_1-3\varpi_2)))
\\
\notag&\simeq
\rH^{i+1}(\cP,\cL_\cP(\begin{tabular}{|c|}
\hline
$\nabla^P(2\varpi_1-2\varpi_2)$
\\
\hline
$-\varpi_2$
\\
\hline
\end{tabular}))=0.
\end{align}

Finally,
\begin{align}
\Ext^i_\cP
&
(\cE(s_2s_1s_2),\cE(w^P))
\simeq
\rH^{i}(\cP,
\cP(\nabla^P(\rho)\otimes
(-2\varpi_2))
\simeq
\rH^{i}(\cP,
\cP(\nabla^P(\varpi_1-\varpi_2))
\\
\notag=0.
\end{align}

\setcounter{equation}{0}
\noindent
(A.6)
Let
$i\in\bbN$.
We have
\begin{align}
\Ext^i_\cP
&
(\cE(s_1s_2s_1s_2),\cE(e))
\simeq
\rH^i(\cP,\cL_\cP(\begin{tabular}{|c|}
\hline
$\nabla^P(2\varpi_1-3\varpi_2)^*$
\\
\hline
$\nabla^P(\varpi_1-2\varpi_2)^*$
\\
\hline
\end{tabular}))
\\
\notag&\simeq
\rH^i(\cP,\cL_\cP(\begin{tabular}{|c|}
\hline
$\nabla^P(2\varpi_1+\varpi_2)$
\\
\hline
$\nabla^P(\rho)$
\\
\hline
\end{tabular}))
\simeq
\delta_{i,0}
\begin{tabular}{|c|}
\hline
$\nabla(2\varpi_1+\varpi_2)$
\\
\hline
$\nabla(\rho)$
\\
\hline
\end{tabular},
\\
\Ext^i_\cP
&
(\cE(s_1s_2s_1s_2),\cE(s_2))
\simeq
\rH^i(\cP,\cL_\cP(\begin{tabular}{|c|}
\hline
$\nabla^P(2\varpi_1+\varpi_2)$
\\
\hline
$\nabla^P(\rho)$
\\
\hline
\end{tabular}
\otimes(-\varpi_2)))
\\
\notag&
\simeq
\rH^i(\cP,\cL_\cP(\begin{tabular}{|c|}
\hline
$\nabla^P(2\varpi_1)$
\\
\hline
$\nabla^P(\varpi_1)$
\\
\hline
\end{tabular}
))
\simeq
\delta_{i,0}
\{\nabla(2\varpi_1)\oplus
\nabla(\varpi_1)\},
\\
\Ext^i_\cP
&
(\cE(s_1s_2s_1s_2),\cE(s_2s_1s_2))
\simeq
\rH^i(\cP,\cL_\cP(\begin{tabular}{|c|}
\hline
$\nabla^P(2\varpi_1+\varpi_2)$
\\
\hline
$\nabla^P(\rho)$
\\
\hline
\end{tabular}
\otimes\nabla^P(\varpi_1-2\varpi_2)))
\\
\notag&
\simeq
\rH^i(\cP,\cL_\cP(\begin{tabular}{|c|}
\hline
$\nabla^P(3\varpi_1-\varpi_2)$
\\
\hline
$\nabla^P(\varpi_1)$
\\
\hline
$\nabla^P(2\varpi_1-\varpi_2)$
\\
\hline
$\varepsilon$
\\
\hline
\end{tabular}
))
\simeq
\delta_{i,0}
\{\nabla(\varpi_1)\oplus
L(e)\}.
\end{align}

One has
\[
\Ext^i_\cP(\cE(s_1s_2s_1s_2),
\cE(s_1s_2))
\simeq
\Ext^i_\cP(\cE(s_1s_2s_1s_2),
\cL_\cP(\begin{tabular}{|c|}
\hline
$\nabla^P(2\varpi_1-2\varpi_2)$
\\
\hline
$\nabla^P(\varpi_1-2\varpi_2)$
\\
\hline
\end{tabular}))
\]
with
$\Ext^i_\cP(\cE(s_1s_2s_1s_2),
\cL_\cP(\nabla^P(\varpi_1-2\varpi_2)))
\simeq
\delta_{i,0}
\{\nabla(\varpi_1)\oplus
\varepsilon\}$ by (3) while
\begin{align*}
\Ext^i_\cP
&
(\cE(s_1s_2s_1s_2),\cL_\cP(
\nabla^P(2\varpi_1-2\varpi_2))\\
\notag&\simeq
\rH^i(\cP,\cL_\cP(\begin{tabular}{|c|}
\hline
$\nabla^P(2\varpi_1+\varpi_2)$
\\
\hline
$\nabla^P(\rho)$
\\
\hline
\end{tabular}
\otimes\nabla^P(2\varpi_1-2\varpi_2)))
\\
\notag&
\simeq
\rH^i(\cP,\cL_\cP(\begin{tabular}{|c|}
\hline
$\nabla^P(4\varpi_1-\varpi_2)$
\\
\hline
$\nabla^P(2\varpi_1)$
\\
\hline
$\varpi_2$
\\
\hline
$\nabla^P(3\varpi_1-\varpi_2)$
\\
\hline
$\nabla^P(\varpi_1)$
\\
\hline
\end{tabular}
))
\simeq
\delta_{i,0}
\{\nabla(2\varpi_1)
\oplus
\nabla(\varpi_2)
\oplus
\nabla(\varpi_1)
\}.
\end{align*}
It follows that
\begin{align}
\Ext^i_\cP(\cE(s_1s_2s_1s_2),
\cE(s_1s_2))
\simeq
\delta_{i,0}
\{\nabla(2\varpi_1)
\oplus
\nabla(\varpi_2)
\oplus\nabla(\varpi_1)^{\oplus_2}\oplus
L(e)\}.
\end{align}

One has
\begin{align*}
\Ext^i_\cP
&
(\cE(s_1s_2s_1s_2),\cE(s_1s_2s_1s_2))
\simeq
\Ext^{i}_\cP(\cL_\cP(\begin{tabular}{|c|}
\hline
$\nabla^P(\varpi_1-2\varpi_2)$
\\
\hline
$\nabla^P(2\varpi_1-3\varpi_2)$
\\
\hline
\end{tabular}),
\cE(s_1s_2s_1s_2),
))
\\
\notag&\simeq
\Ext^{i}_\cP(\cL_\cP(\nabla^P(2\varpi_1-3\varpi_2)),
\cE(s_1s_2s_1s_2)
)
\quad\text{by (A.5.6)},
\end{align*}
the right hand side of which fits by (A.1.2)
into the long exact sequence (A.5.5)
with
$\cE(s_2s_1s_2)$ replaced by
$\cL_\cP(\nabla^P(2\varpi_1-3\varpi_2))$.
As
$\Ext^\bullet_\cP(
\cL_\cP(
\nabla^P(2\varpi_1-3\varpi_2),
\cL_\cP((-2\varpi_2)
\otimes
\nabla(\varpi_1)))
\simeq
\rH^\bullet(\cP,
\cL_\cP(
\nabla^P(2\varpi_1-\varpi_2)))
\otimes
\nabla(\varpi_1)
=0$, 
one obtains
\begin{align}
\Ext^i_\cP
&
(\cE(s_1s_2s_1s_2),\cE(s_1s_2s_1s_2))
\\
\notag&\simeq
\Ext^{i+1}_\cP(\cL_\cP(\nabla^P(2\varpi_1-3\varpi_2)),
\cL_\cP(\nabla^P(\varpi_1-3\varpi_2)))
\\
\notag&\simeq
\rH^{i+1}(\cP,
\cL_\cP(\nabla^P(2\varpi_1-3\varpi_2)^*\otimes
\nabla^P(\varpi_1-3\varpi_2)))
\\
\notag&\simeq
\rH^{i+1}(\cP,
\cL_\cP(\nabla^P(2\varpi_1+\varpi_2)\otimes
\nabla^P(\varpi_1-3\varpi_2)))
\\
\notag&\simeq
\rH^{i+1}(\cP,
\cL_\cP(\begin{tabular}{|c|}
\hline
$\nabla^P(3\varpi_1-2\varpi_2)$
\\
\hline
$\nabla^P(\varpi_1-\varpi_2)$
\\
\hline
\end{tabular}))
\simeq
\rH^{i+1}(\cP,
\cL_\cP(\nabla^P(3\varpi_1-2\varpi_2)))
\\
\notag&=
\rH^{i+1}(\cP,
\cL_\cP(\nabla^P(s_1\bullet0)))
\simeq
\delta_{i+1,1}\Bbbk
\quad\text{by
\cite[II.5.5]{J}}
\\
\notag&=
\delta_{i,0}\Bbbk.
\end{align}

Finally, one has
\begin{align}
\Ext^i_\cP
&
(\cE(s_1s_2s_1s_2),\cE(w^P))
\simeq
\rH^i(\cP,\cL_\cP(\begin{tabular}{|c|}
\hline
$\nabla^P(2\varpi_1+\varpi_2)$
\\
\hline
$\nabla^P(\rho)$
\\
\hline
\end{tabular}
\otimes(-2\varpi_2)))
\quad\text{as in (1)}
\\
\notag&
\simeq
\rH^i(\cP,\cL_\cP(\begin{tabular}{|c|}
\hline
$\nabla^P(2\varpi_1-\varpi_2)$
\\
\hline
$\nabla^P(\varpi_1-\varpi_2)$
\\
\hline
\end{tabular}
))
=0.
\end{align}

\setcounter{equation}{0}
\noindent
(A.7)
Let
$i\in\bbN$.
We have
\begin{align}
\Ext^i_\cP
&
(\cE(w^P),\cE(e))
\simeq
\rH^i(\cP, \cL_\cP(2\varpi_2))
\simeq
\delta_{i,0}\nabla(2\varpi_2),
\\
\Ext^i_\cP
&
(\cE(w^P),\cE(s_2))
\simeq
\rH^i(\cP, \cL_\cP(\varpi_2))
\simeq
\delta_{i,0}\nabla(\varpi_2),
\\
\Ext^i_\cP
&
(\cE(w^P),\cE(s_1s_2))
\simeq
\rH^i(\cP, \cL_\cP(\begin{tabular}{|c|}
\hline
$\nabla^P(2\varpi_1)$
\\
\hline
$\nabla^P(\varpi_1)$
\\
\hline
\end{tabular}))
\\
\notag&\simeq
\delta_{i,0}\{\nabla(2\varpi_1)\oplus\nabla(\varpi_1)\},
\\
\Ext^i_\cP
&
(\cE(w^P),\cE(s_2s_1s_2))
\simeq
\rH^i(\cP, \cL_\cP(\nabla^P(\varpi_1)))
\simeq
\delta_{i,0}\nabla(\varpi_1),
\\
\Ext^i_\cP
&
(\cE(w^P),\cE(s_1s_2s_1s_2))
\simeq
\rH^i(\cP, \cL_\cP(\begin{tabular}{|c|}
\hline
$\nabla^P(\varpi_1)$
\\
\hline
$\nabla^P(2\varpi_1-\varpi_2)$
\\
\hline
\end{tabular}))
\simeq
\delta_{i,0}\nabla(\varpi_1),
\\
\Ext^i_\cP
&
(\cE(w^P),\cE(w^P))
\simeq
\rH^i(\cP, \cL_\cP(\varepsilon))
\simeq
\delta_{i,0}\Bbbk.
\end{align}

\setcounter{equation}{0}
\noindent
(A.8)
We have thus shown 
\begin{prop}
Assume $p\geq11$.
The $\cE(w)$'s, $w\in W^P$, 
form a strongly exceptional sequence on $\cP$ such that
$\Mod_\cP(\cE(x),\cE(y))\ne0$ iff  $x\geq y$ in the Chevalley-Bruhat order with
isomorphisms of $G$-modules
\begin{align*}
\Mod_\cP
(\cE(s_2),\cE(e))
&\simeq
\nabla(\varpi_2),
\quad
\Mod_\cP
(\cE(s_1s_2),\cE(e))
\simeq
\nabla(\rho)\oplus\nabla(2\varpi_1),
\\
\Mod_\cP
(\cE(s_1s_2),\cE(s_2))
&\simeq
\nabla(\varpi_1),
\quad
\Mod_\cP
(\cE(s_2s_1s_2),\cE(e))
\simeq
\nabla(\rho),
\\
\Mod_\cP
(\cE(s_2s_1s_2),\cE(s_2))
&\simeq
\nabla(\varpi_1),
\quad
\Mod_\cP
(\cE(s_2s_1s_2),\cE(s_1s_2))
\simeq
\nabla(\varpi_1)\oplus
L(e),
\\
\Mod_\cP
(\cE(s_1s_2s_1s_2),\cE(e))
&\simeq
\begin{tabular}{|c|}
\hline
$\nabla(2\varpi_1+\varpi_2)$
\\
\hline
$\nabla(\rho)$
\\
\hline
\end{tabular},
\\
\notag
\Mod_\cP
(\cE(s_1s_2s_1s_2),\cE(s_2))
&\simeq
\nabla(2\varpi_1)\oplus
\nabla(\varpi_1),
\\
\notag
\Mod_\cP
(\cE(s_1s_2s_1s_2),\cE(s_1s_2))
&\simeq
\nabla(2\varpi_1)
\oplus
\nabla(\varpi_2)
\oplus\nabla(\varpi_1)^{\oplus_2}\oplus
L(e),
\\
\notag
\Mod_\cP
(\cE(s_1s_2s_1s_2),\cE(s_2s_1s_2))
&\simeq
\nabla(\varpi_1)
\oplus
L(e),
\\
\Mod_\cP
(\cE(w^P),\cE(e))
&\simeq
\nabla(2\varpi_2),
\quad
\Mod_\cP
(\cE(w^P),\cE(s_2))
\simeq
\nabla(\varpi_2),
\\
\Mod_\cP
(\cE(w^P),\cE(s_1s_2))
&\simeq
\nabla(2\varpi_1)\oplus\nabla(\varpi_1),
\\
\Mod_\cP
(\cE(w^P),\cE(s_2s_1s_2))
&\simeq
\nabla(\varpi_1)
\simeq
\Mod_\cP
(\cE(w^P),\cE(s_1s_2s_1s_2)).
\end{align*}

\end{prop}

\setcounter{equation}{0}
\noindent
(A.9)
To see that the $\cE(w)$,
$w\in
W^P$,
Karoubian generate $\rD^b(\coh\cP)$,
it is enough by a result
attributed to Kontsevich
by Positselskii
\cite[Th. 3.5.1]{BMR02}
to verify that all
$\cL_\cP(-2n\rho_P)$, $n\in\bbN^+$,
are Karoubian generated by the 
$\cE(w)$'s,
where
$2\rho_P=\sum_{\beta\in R^+\setminus\{\alpha_1\}}\beta=
2\varpi_2$
in the present setting.

Let
$\hat\cE=\langle\cE(w)
\mid
w\in
W^P\rangle$ denote the triangulated subcategory of
$\rD^b(\coh\cP)$
Karoubian generated by
$\cE(w)$,
$w\in
W^P$.
Let
$\pi:\cB\to\cP$ be the natural morphism.
Recall from \cite[1.3.6]{Or}
the projection formula
$\id_{\rD^b(\coh\cP)}\simeq
(\bfR\pi_*)\circ\pi^*:\rD^b(\coh\cP)\to\rD^b(\coh\cP)$.
Recall also from
\cite[I.5.17]{J} an isomorphism
$\pi^*\cL_\cP(M)\simeq\cL(M)$ for each $P$-module $M$, and from
\cite[I.5.19]{J}
an isomorphism
$\bfR\pi_*\cL(m\varpi_1+n\varpi_2)\simeq
\cL_\cP(\nabla^P(m\varpi_1+n\varpi_2))$
for each $m\in\bbN$ and $n\in\bbZ$.
Then, setting 
$\tilde\cE=\langle\pi^*\cE(w)|w\in W^P\rangle$ to be the 
triangulated subcategory of
$\rD^b(\coh\cB)$
Karoubian generated by
$\pi^*\cE(w)$,
$w\in
W^P$, it is enough to show that
\begin{align}
\cL(-m\varpi_2)\in\tilde\cE
\quad\forall m\in\bbN.\end{align}
For our purpose
note also that
whenever
$\cL(\nabla^P(m\varpi_1+n\varpi_2))\in\tilde\cE$, $m\in\bbN$, $n\in\bbZ$,
we may also allow
$\cL(m\varpi_1+n\varpi_2)\in\tilde\cE$, $m\in\bbN$, and conversely.

We thus start with
$\tilde\cE=\langle
\cO_\cB,
\cL(-\varpi_2),
\cL(\nabla^P(\varpi_1-2\varpi_2)),
\cL(\varpi_1-2\varpi_2),
\cL(\nabla^P(2\varpi_1-2\varpi_2)),
\cL(2\varpi_1-2\varpi_2),
\cL(\nabla^P(2\varpi_1-3\varpi_2)),
\cL(2\varpi_1-3\varpi_2),
\cL(-2\varpi_2)
\rangle$.
In view of the $B$-filtration on
$\nabla^{P}(\varpi_1-2\varpi_2)=
\begin{tabular}{|c|}
\hline
$\varpi_1-2\varpi_2$
\\
\hline
$-\rho$
\\
\hline
\end{tabular}$,
as
$\cL(\nabla^{P}(\varpi_1-2\varpi_2))$
and
$\cL(\varpi_1-2\varpi_2)\in\tilde\cE$, one has also
$\cL(-\rho)\in\tilde\cE$.
Likewise
from the $B$-filtration
$\nabla^{P}(2\varpi_1-2\varpi_2)=
\begin{tabular}{|c|}
\hline
$2\varpi_1-2\varpi_2$
\\
\hline
$-\varpi_2$
\\
\hline
$-2\varpi_1$
\\
\hline
\end{tabular}$,
get
$\cL(-2\varpi_1)\in\tilde\cE$.
Then from
$\nabla^{P}(2\varpi_1-3\varpi_2)=
\begin{tabular}{|c|}
\hline
$2\varpi_1-3\varpi_2$
\\
\hline
$-2\varpi_2$
\\
\hline
$-2\varpi_1-\varpi_2$
\\
\hline
\end{tabular}$,
get
$\cL(-2\varpi_1-\varpi_2)\in\tilde\cE$.

We now explain our strategy.
Recall from
\cite[5.1.4]{HKR} the Koszul resolution
\begin{multline}
0\to
\cL(-7\varpi_1)\otimes\wedge^7\nabla(\varpi_1)
\to
\cL(-6\varpi_1)\otimes\wedge^6\nabla(\varpi_1)
\to
\dots
\\
\to
\cL(-\varpi_1)\otimes\wedge^1\nabla(\varpi_1)
\to
\cO_\cB
\to
0.
\end{multline}
Tensoring entries of
$\tilde\cE$ with $\nabla(\varpi_1)$ and
$\nabla(\varpi_2)$, we will find
$\cL(r\varpi_1-\varpi_2)\in\tilde\cE$, $r\in[-3,3]$.
Then from (2)
we obtain all
$\cL(n\varpi_1-\varpi_2)\in\tilde\cE$, 
$n\in\bbZ$, and hence also
all
$\cL(\nabla^P(m\varpi_1-\varpi_2))\in\tilde\cE$, 
$m\in\bbN$.
That, in turn, will yield all
$\cL(r\varpi_1)\in\tilde\cE$, $r\in[-4,2]$,
and hence by (2) again all
$\cL(n\varpi_1)$, 
$n\in\bbZ$, 
in
$\tilde\cE$.
Then get all
$\cL(r\varpi_1-2\varpi_2)\in\tilde\cE$, $r\in[-2,4]$,
and hence all
$\cL(n\varpi_1-2\varpi_2)\in\tilde\cE$,
$n\in\bbZ$.
Thus,
all the weights $\nu$
of $\nabla(\varpi_2)\otimes
(-\varpi_2)$
will have
$\cL(\nu)\in\tilde\cE$,
and hence
tensoring $\cL(n\varpi_1-\varpi_2)$ with
$\nabla(\varpi_2)\otimes
\cL(-\varpi_2)$ will yield all
$\cL(n\varpi_1-3\varpi_2))\in\tilde\cE$,
$n\in\bbZ$.
Repeat the procedure 
to obtain all 
$\cL(n\varpi_1-m\varpi_2)\in\tilde\cE$,
$n\in\bbZ$, $m\in\bbN$,
and hence in particular
(1).

\setcounter{equation}{0}
\noindent
(A.10)
In
(A.9)
we have added 
$\cL(-\rho), \cL(-2\varpi_1)$ and
$\cL(-2\varpi_1-\varpi_2)$
to $\tilde\cE$, so that
$\tilde\cE=\langle
\cO_\cB,
\cL(-\varpi_2),
\cL(\nabla^P(\varpi_1-2\varpi_2)),
\cL(\varpi_1-2\varpi_2),
\cL(\nabla^P(2\varpi_1-2\varpi_2)),
\cL(2\varpi_1-2\varpi_2),
\cL(\nabla^P(2\varpi_1-3\varpi_2)),
\cL(2\varpi_1-3\varpi_2),
\cL(-2\varpi_2),
\cL(-\rho),
\cL(-2\varpi_1),
\cL(-2\varpi_1-\varpi_2)
\rangle$.

As $\cL(-\varpi_2)\in\tilde\cE$,
so does
\begin{align*}
\nabla(\varpi_1)\otimes_\Bbbk
\cL(-\varpi_2)
&\simeq
\cL(\nabla(\varpi_1)\otimes
(-\varpi_2))
\simeq
\cL(\begin{tabular}{|c|}
\hline
$\nabla^{P}(\varpi_1-\varpi_2)$
\\
\hline
$\nabla^{P}(2\varpi_1-2\varpi_2)$
\\
\hline
$\nabla^{\alpha_1}(\varpi_1-2\varpi_2)$
\\
\hline
\end{tabular}).
\end{align*}
As
$\cL(\nabla^{P}(\varpi_1-2\varpi_2))$
and $\cL(\nabla^{P}(2\varpi_1-2\varpi_2))
\in\tilde\cE$,
one has
$\cL(\nabla^{P}(\varpi_1-\varpi_2))\in\tilde\cE$, and hence also
$\cL(\varpi_1-\varpi_2)\in\tilde\cE$.
Then, from
$\nabla^{P}(\varpi_1-\varpi_2)
=
\begin{tabular}{|c|}
\hline
$\varpi_1-\varpi_2$
\\
\hline
$-\varpi_1$
\\
\hline
\end{tabular}$, get
also
$\cL(-\varpi_1)\in\tilde\cE$.

Likewise, as
$\cL(-2\varpi_2)\in\tilde\cE$,
get from
$\nabla(\varpi_1)\otimes(-2\varpi_2)\simeq
\begin{tabular}{|c|}
\hline
$\nabla^{P}(\varpi_1-2\varpi_2)$
\\
\hline
$\nabla^{P}(2\varpi_1-3\varpi_2)$
\\
\hline
$\nabla^{P}(\varpi_1-3\varpi_2)$
\\
\hline
\end{tabular})$
that
$\cL(\nabla^{P}(\varpi_1-3\varpi_2))$
and
$\cL(\varpi_1-3\varpi_2)\in\tilde\cE$.
Then, from
$\nabla^{P}(\varpi_1-3\varpi_2)
=
\begin{tabular}{|c|}
\hline
$\varpi_1-3\varpi_2$
\\
\hline
$-\varpi_1-2\varpi_2$
\\
\hline
\end{tabular}$, get
also
$\cL(-\varpi_1-2\varpi_2)\in\tilde\cE$.

One has
$\nabla(\varpi_1)\otimes\cL(\varpi_1-2\varpi_2)\in\tilde\cE$.
As all weights
$\nu$ of
$\nabla(\varpi_1)\otimes(\varpi_1-2\varpi_2)
$ except $3\varpi_1-3\varpi_2$
have
$\cL(\nu)\in\tilde\cE$,
$\cL(3\varpi_1-3\varpi_2)\in\tilde\cE$, and hence also
$\cL(\nabla^P(3\varpi_1-3\varpi_2))\in\tilde\cE$.
Then,
as
$\nabla^P(3\varpi_1-3\varpi_2)
=
\begin{tabular}{|c|}
\hline
$3\varpi_1-3\varpi_2$
\\
\hline
$\varpi_1-2\varpi_2$
\\
\hline
$-\varpi_1-\varpi_2$
\\
\hline
$-3\varpi_1$
\\
\hline
\end{tabular}$,
get also
$\cL(-3\varpi_1)\in\tilde\cE$.

Now,
$\nabla(\varpi_2)\otimes\cL(-\varpi_2)\in\tilde\cE$
with
\[
\nabla(\varpi_2)\otimes(-\varpi_2)
=
\begin{tabular}{|c|}
\hline
$\varepsilon$
\\
\hline
$\nabla^{P}(3\varpi_1-2\varpi_2)$
\\
\hline
$\nabla^{\alpha_1}(2\varpi_1-2\varpi_2)$
\\
\hline
$-\varpi_2$
\\
\hline
$\nabla^{\alpha_1}(3\varpi_1-3\varpi_2)$
\\
\hline
$-2\varpi_2$
\\
\hline
\end{tabular}.
\]
As all components of
$\nabla(\varpi_2)\otimes\cL(-\varpi_2)$ except $\cL(
\nabla^{P}(3\varpi_1-2\varpi_2))$
belong to
$\tilde\cE$,
so does
$\cL(
\nabla^{P}(3\varpi_1-2\varpi_2))$,
and hence also
$\cL(3\varpi_1-2\varpi_2)$.
As
$\nabla^{P}(3\varpi_1-2\varpi_2)=
\begin{tabular}{|c|}
\hline
$3\varpi_1-2\varpi_2$
\\
\hline
$\varpi_1-\varpi_2$
\\
\hline
$-\varpi_1$
\\
\hline
$-3\varpi_1+\varpi_2$
\\
\hline
\end{tabular}$
has all its weights 
$\nu$
but
$-3\varpi_1+\varpi_2$
such that
$\cL(\nu)\in\tilde\cE$,
$\cL(-3\varpi_1+\varpi_2)\in\tilde\cE$.

We have seen above that
$\cL(\varpi_1-\varpi_2)\in\tilde\cE$, and hence also
$\nabla(\varpi_1)\otimes
\cL(\varpi_1-\varpi_2)\in\tilde\cE$.
As all weights
$\nu$ of
$\nabla(\varpi_1)\otimes(\varpi_1-\varpi_2)$ except $2\varpi_1-\varpi_2$
have
$\cL(\nu)\in\tilde\cE$,
$\cL(2\varpi_1-\varpi_2)\in\tilde\cE$,
and hence also
$\cL(\nabla^P(2\varpi_1-\varpi_2))\in\tilde\cE$.
Then,
as
$\nabla^P(2\varpi_1-\varpi_2)=
\begin{tabular}{|c|}
\hline
$2\varpi_1-\varpi_2$
\\
\hline
$\varepsilon$
\\
\hline
$-2\varpi_1+\varpi_2$
\\
\hline
\end{tabular}$,
get also
$\cL(-2\varpi_1+\varpi_2)\in\tilde\cE$.

One has
$\nabla(\varpi_1)\otimes
\cO_\cB\in\tilde\cE$
with
$
\nabla(\varpi_1)=
\begin{tabular}{|c|}
\hline
$\nabla^P(\varpi_1)$
\\
\hline
$\nabla^P(2\varpi_1-\varpi_2)$
\\
\hline
$\nabla^P(\varpi_1-\varpi_2)$
\\
\hline
\end{tabular}$.
As all components except
$\cL(\nabla^P(\varpi_1))$ belong to
$\tilde\cE$, so does
$\cL(\nabla^P(\varpi_1))$, and hence also
$\cL(\varpi_1)$.
Then, as
$
\nabla^P(\varpi_1)=
\begin{tabular}{|c|}
\hline
$\varpi_1$
\\
\hline
$-\varpi_1+\varpi_2$
\\
\hline
\end{tabular}$,
$\cL(-\varpi_1+\varpi_2)\in\tilde\cE$ also.

As
$\nabla(\varpi_1)\otimes\cL(2\varpi_1-2\varpi_2)\in\tilde\cE$
and as all weights
$\nu$ of
$\nabla(\varpi_1)\otimes(2\varpi_1-2\varpi_2)
$ except $4\varpi_1-3\varpi_2$
have
$\cL(\nu)\in\tilde\cE$,
$\cL(4\varpi_1-3\varpi_2)\in\tilde\cE$, and hence also
$\cL(\nabla^P(4\varpi_1-3\varpi_2))\in\tilde\cE$.
Then, 
as
all weights $\nu$ of
$\nabla^P(4\varpi_1-3\varpi_2)$ except 
$-4\varpi_1+\varpi_2$
have $\cL(\nu)\in\tilde\cE$,
$\cL(-4\varpi_1+\varpi_2)\in\tilde\cE$.

Now,
$\nabla(\varpi_2)\otimes\cL(\varpi_1-\varpi_2)\in\tilde\cE$.
As all weights $\nu$
of
$\nabla(\varpi_2)\otimes(\varpi_1-\varpi_2)$ except
$4\varpi_1-2\varpi_2$
have
$\cL(\nu)\in\tilde\cE$,
$\cL(4\varpi_1-2\varpi_2)\in\tilde\cE$.
Likewise,
as $\cL(-\rho)\in\tilde\cE$,
so does
$\nabla(\varpi_2)\otimes\cL(-\rho)$.
As all weights $\nu$
of
$\nabla(\varpi_2)\otimes(-\rho)$ except
$-4\varpi_1$
have
$\cL(\nu)\in\tilde\cE$, get
$\cL(-4\varpi_1)\in\tilde\cE$.

One has
$\nabla(\varpi_1)\otimes\cL(-2\varpi_1-\varpi_2)\in\tilde\cE$.
As all weights $\nu$ of
$\nabla(\varpi_1)\otimes(-2\varpi_1-\varpi_2)$ except 
$-3\varpi_1-\varpi_2$ have
$\cL(\nu)\in\tilde\cE$, 
so does
$\cL(-3\varpi_1-\varpi_2)$.
Likewise, we have seen
$\cL(2\varpi_1-\varpi_2)\in\tilde\cE$, and hence
$\nabla(\varpi_1)\otimes\cL(2\varpi_1-\varpi_2)\in\tilde\cE$.
As all weights $\nu$ of
$\nabla(\varpi_1)\otimes(2\varpi_1-\varpi_2)$ except 
$3\varpi_1-\varpi_2$ have
$\cL(\nu)\in\tilde\cE$, 
so does
$\cL(3\varpi_1-\varpi_2)$.
At this point we have
all $\cL(r\varpi_1-\varpi_2)\in\tilde\cE$,
$r\in[-3,3]$.
Then, using (A.9.2),
one obtains 
\begin{align}
\cL(n\varpi_1-\varpi_2)\in\tilde\cE
\quad\forall n\in\bbZ, \text{ and hence also }
\cL(\nabla^P(m\varpi_1-\varpi_2))\in\tilde\cE
\quad\forall m\in\bbZ.
\end{align}

Then, 
as 
$\nabla(\varpi_2)
=
\begin{tabular}{|c|}
\hline
$\varpi_2$
\\
\hline
$\nabla^P(3\varpi_1-\varpi_2)$
\\
\hline
$\nabla^P(2\varpi_1-\varpi_2)$
\\
\hline
$\varepsilon$
\\
\hline
$\nabla^P(3\varpi_1-2\varpi_2)$
\\
\hline
$-\varpi_2$
\\
\hline
\end{tabular}$,
all components of
$\nabla(\varpi_2)\otimes\cO_\cB\in\tilde\cE$
except $\cL(\varpi_2)$
belong to $\tilde\cE$, so does
$\cL(\varpi_2)$.
In turn,
$\nabla(\varpi_1)\otimes\cL(\varpi_1)\in\tilde\cE$.
As all weights $\nu$
of
$\nabla(\varpi_1)\otimes\varpi_1$
except $2\varpi_1$
have
$\cL(\nu)\in\tilde\cE$,
one has
$\cL(2\varpi_1)\in\tilde\cE$.
Thus,
all
$\cL(r\varpi_1)\in\tilde\cE$, $r\in[-4,2]$.
Then, using (A.9.2) again, get
\begin{align}
\cL(n\varpi_1)\in\tilde\cE
\quad\forall n\in\bbZ. 
\end{align}

By (1) one has
$\nabla(\varpi_1)\otimes\cL(-3\varpi_1-\varpi_2)\in\tilde\cE$.
As all weights $\nu$
of
$\nabla(\varpi_1)\otimes
(-3\varpi_1-\varpi_2)$
except $-2\varpi_1-2\varpi_2$
have
$\cL(\nu)\in\tilde\cE$,
one has
$\cL(-2\varpi_1-2\varpi_2)\in\tilde\cE$.
Then
all
$\cL(r\varpi_1-2\varpi_2)\in\tilde\cE$, $r\in[-2,4]$.
It follows from (A.9.2) that
\begin{align}
\cL(n\varpi_1-2\varpi_2)\in\tilde\cE
\quad\forall n\in\bbZ. 
\end{align}

Then,
$\nabla(\varpi_1)\otimes\cL(-\varpi_1-2\varpi_2)\in\tilde\cE$.
As all weights $\nu$
of
$\nabla(\varpi_1)\otimes
(-\varpi_1-2\varpi_2)$
except $-3\varpi_2$
have
$\cL(\nu)\in\tilde\cE$ by (1) and (3),
one has
$\cL(-3\varpi_2)\in\tilde\cE$.
Likewise,
$\nabla(\varpi_1)\otimes\cL(n\varpi_1-2\varpi_2)\in\tilde\cE$,
$n\in\bbZ$, will yield
\begin{align}
\cL(n\varpi_1-3\varpi_2)\in\tilde\cE
\quad\forall n\in\bbZ.
\end{align}

In order to get
all
$\cL(n\varpi_1-4\varpi_2)\in\tilde\cE$, $n\in\bbZ$,
consider
$\nabla(\varpi_2)\otimes
\cL(n\varpi_1-2\varpi_2)\in\tilde\cE$
by (3).
All the weights $\nu$
of
$\nabla(\varpi_2)\otimes
(n\varpi_1-2\varpi_2)$ except
$(n+3)\varpi_1-4\varpi_2$
have 
$\cL(\nu)\in\tilde\cE$ by
(1), (3) and (4), and hence also
$\cL((n+3)\varpi_1-4\varpi_2)\in\tilde\cE$.
To see
all
$\cL(n\varpi_1-5\varpi_2)\in\tilde\cE$,
use
$\nabla(\varpi_2)\otimes
\cL(n\varpi_1-3\varpi_2)\in\tilde\cE$ from (4) to obtain
$\cL((n+3)\varpi_1-5\varpi_2)\in\tilde\cE$.
Repeat to get all
$\cL(\nabla^P(n\varpi_1-m\varpi_2))\in\tilde\cE$, $m\in\bbN$, $n\in\bbZ$, as desired.

\setcounter{equation}{0}
\noindent
(A.10)
We have thus obtained

\begin{thm}
Assume $p\geq11$.
The $\cE(w)$, $w\in W^P$,
form a Karoubian complete strongly exceptional collection on $\cP$ such that
$\Mod_\cP(\cE(x),\cE(y))\ne0$
iff
$x\geq y$ in the Chevalley-Bruhat order.

\end{thm}

\noindent
{\bf B.}
Let us also write down an easier case of the parabolic $P$ associated to the long simple root; parametrization of the sheaves is, as in {\bf A}, different from the one in
\cite{KY} twisted by $w_0?w_P$.
We assume that
$p\geq11$.

\setcounter{equation}{0}
\noindent
(B.1)
Let 
$
\soc^i\hat\nabla_P(\varepsilon)$
denote the $i$-th $G_1T$-socle of $\hat\nabla_P(\varepsilon)=\ind_{G_1P}^G(\varepsilon)$, and put
$\soc_i\hat\nabla_P(\varepsilon)=
\soc^i\hat\nabla_P(\varepsilon)/\soc^{i-1}\hat\nabla_P(\varepsilon)$.
From its $G_1T$-structure we readily obtain
\begin{align}
\soc_1\hat\nabla_P(\varepsilon)&\simeq
L(e)\otimes\varepsilon,
\\
\notag
\soc_2\hat\nabla_P(\varepsilon)&\simeq
L(s_1)\otimes
(-\varpi_1)^{[1]},
\\
\notag
\soc_3\hat\nabla_P(\varepsilon)&\simeq
L(s_2s_1)\otimes
(-2\varpi_1)^{[1]}
\oplus
L(e)\otimes(-2\varpi_1)^{[1]},
\\
\notag
\soc_4\hat\nabla_P(\varepsilon)&\simeq
L(s_1s_2s_1)\otimes
\{(-3\varpi_1)\otimes
\begin{tabular}{|c|}
\hline
$\varpi_1$
\\
\hline
$\nabla^P(-\varpi_1+\varpi_2)$
\\
\hline
$\varepsilon$
\\
\hline
\end{tabular}
\}^{[1]}
\\
\notag&
\oplus
L(w^P)\otimes(-2\varpi_1)^{[1]}
\oplus
L(s_1)\otimes(-3\varpi_1)^{[1]},
\\
\notag
\soc_5\hat\nabla_P(\varepsilon)&\simeq
L(s_2s_1s_2s_1)\otimes
(-3\varpi_1)^{[1]}\oplus
L(e)\otimes(-4\varpi_1)^{[1]},
\\
\notag
\soc_6\hat\nabla_P(\varepsilon)&\simeq
L(
w^P)\otimes
(-4\varpi_1)^{[1]},
\end{align}
where $\nabla^P=\ind_B^P$ and $w^P=s_1s_2s_1s_2s_1$.

For
$w\in W^P$ we let $L(w)\otimes(\soc^1_{i,w})^{[1]}$
denote the 
$L(w)$-isotypic part of $\soc_{i}\hat\nabla_P(\varepsilon)$.
We show that
$3\varpi_1\otimes\soc^1_{4,s_1s_2s_1}
=
\begin{tabular}{|c|}
\hline
$\varpi_1$
\\
\hline
$\nabla^P(-\varpi_1+\varpi_2)$
\\
\hline
$\varepsilon$
\\
\hline
\end{tabular}
$
is an indecomposable $P$-module isomorphic to
the quotient of $\nabla(\varpi_1)$
by a $P$-submodule generated by 
a vector of weight
$-2\varpi_1+\varpi_2$.

\setcounter{equation}{0}
\noindent
(B.2)
Just suppose 
$(-3\varpi_1)\otimes
\nabla^P(-\varpi_1+\varpi_2)
\simeq
\nabla^P(-4\varpi_1+\varpi_2)$ is a $P$-submodule of 
$\soc^1_{4,s_1s_2s_1}$.
Then,
$-3\varpi_1$ would be a direct summand of
$\soc^1_{4,s_1s_2s_1}$
as
$\Ext^1_P(-2\varpi_1,-3\varpi_1)=0$.
Dualizing, $L(s_1s_2s_1)\otimes
(3\varpi_1)$ would be a direct summand of
the third subquotient 
$\rad_3(\hat\nabla_P(\varepsilon)^*)=
\rad^3(\hat\nabla_P(\varepsilon)^*)/\rad^4(\hat\nabla_P(\varepsilon)^*)$
in the $G_1T$ radical series of
$\hat\nabla_P(\varepsilon)^*\simeq
\hat\nabla_P(5(p-1)\varpi_1)$.
By \cite{AbK} one has 
$\rad^i\hat\nabla_P(5(p-1)\varpi_1)
=
\soc^{6-i}\hat\nabla_P(5(p-1)\varpi_1)$.
It would follow that there is a $P$-submodule $M$ of
$\soc^3\hat\nabla_P(5(p-1)\varpi_1)$ containing
$\soc^2\hat\nabla_P(5(p-1)\varpi_1)$ such that
$M/\soc^2\hat\nabla_P(5(p-1)\varpi_1)\simeq
L(s_1s_2s_1)\otimes3p\varpi_1$.
That would induce an exact sequence of
$G$-modules
\[
\ind_{G_1P}^G(M\otimes(-3p\varpi_1))
\to
L(s_1s_2s_1)\to
\rR^1\ind_{G_1P}^G(\soc^2(\hat\nabla_P((2p-5)\varpi_1))).
\]
But the $G_1P$-components of
$\soc^2\hat\nabla_P((2p-5)\varpi_1)$ 
are just $L(w^P)\otimes
p\varpi_1$,
$L(e)\otimes
p\varpi_1$, and $L(s_2s_1s_2s_1)$
\cite[1.6.12]{KY}, and hence
$\rR^1\ind_{G_1P}^G(\soc^2(\hat\nabla_P((2p-5)\varpi_1)))=0$.
Also,
\begin{align*}
\ind_{G_1P}^G(M\otimes(-3p\varpi_1))&\leq
\ind_{G_1P}^G(\hat\nabla_P(5(p-1)\varpi_1)\otimes(-3p\varpi_1))
\simeq
\nabla((2p-5)\varpi_1)
\end{align*}
with
$\nabla((2p-5)\varpi_1)
$ not having a $G$-composition factor
$L(s_1s_2s_1)$
\cite[p. 149]{A86}, absurd.

Now, 
\begin{align*}
\Ext^1_P
&
(\nabla^P(-\varpi_1+\varpi_2),\varepsilon)
\simeq
\rH^1(P,\nabla^P(-\varpi_1+\varpi_2)^*)
\simeq
\rH^1(P,\nabla^P(-2\varpi_1+\varpi_2))
\\
&\simeq
\Ext^1_P(\varepsilon,\nabla^P(-2\varpi_1+\varpi_2))
\simeq
\Ext^1_B(\varepsilon,-2\varpi_1+\varpi_2)
\\
&\simeq
\Mod_G(\varepsilon,\rR^1\ind_B^G(-2\varpi_1+\varpi_2)
)
=
\Mod_G(\varepsilon,\rR^1\ind_B^G(s_1\bullet0)
)
\simeq
\nabla(\varepsilon)=\Bbbk.
\end{align*}
Thus,
$3\varpi_1\otimes\soc^1_{4,s_1s_2s_1}$ has, up to isomorphism, a unique indecomposable $P$-submodule $E$
extending
$\nabla^P(-\varpi_1+\varpi_2)$ by
$\varepsilon$.

Next,
just suppose
$-2\varpi_1\leq
\soc^1_{4,s_1s_2s_1}$.
Then there would be a $P$-submodule $M'$ of 
$\soc^4\hat\nabla_P(\varepsilon)$
containing
$\soc^3\hat\nabla_P(\varepsilon)$ such that
$M'/\soc^3\hat\nabla_P(\varepsilon)
\simeq
L(s_1s_2s_1)\otimes(-2\varpi_1)^{[1]}$, which would induce an exact sequence of $G$-modules
\[
\ind_{G_1P}^G(M'\otimes2p\varpi_1)
\to
L(s_1s_2s_1)\to
\rR^1\ind_{G_1P}^G(\soc^3\hat\nabla_P(\varepsilon)\otimes2p\varpi_1).
\]
But
$\ind_{G_1P}^G(M'\otimes2p\varpi_1)\leq
\ind_{G_1P}^G(\hat\nabla_P(\varepsilon)\otimes2p\varpi_1)
\simeq
\nabla(2p\varpi_1)$ with
$\nabla(2p\varpi_1)$ having no composition factor
$L(s_1s_2s_1)$
\cite[p. 149]{A86}, and neither
$\rR^1\ind_{G_1P}^G(\soc^3\hat\nabla_P(\varepsilon)\otimes2p\varpi_1)$ has $G_1$-composition factor
$L(s_1s_2s_1)$, absurd.
Thus,
$3\varpi_1\otimes\soc^1_{4,s_1s_2s_1}$ is $P$-indecomposable.
Finally, 
\begin{align*}
\Ext^1_P
&
(\varpi_1,E)
\simeq
\Ext^1_P
(\varpi_1,\nabla^P(-\varpi_1+\varpi_2))
\simeq
\Ext^1_P
(\varepsilon,
\nabla^P(-2\varpi_1+\varpi_2))\simeq
\Bbbk
\end{align*}
as above, and hence
$3\varpi_1\otimes\soc^1_{4,s_1s_2s_1}$ is, up to isomorphism, a unique $P$-extension of
$\varpi_1$ by
$E$.
It follows that
$3\varpi_1\otimes\soc^1_{4,s_1s_2s_1}
\simeq
\nabla(\varpi_1)/
\begin{tabular}{|c|}
\hline
$\nabla^P(-2\varpi_1+\varpi_2)$
\\
\hline
$-\varpi_1$
\\
\hline
\end{tabular}
$.

\begin{prop}
All $\soc^1_{\ell(w)+1,w}$, $w\in W^P$,  
are
indecomposable $P$-modules,
of highest weight
$w^{-1}\bullet(w\bullet0)^1$
except for $w=s_1s_2s_1$.
In the last case it is generated by a vector of weight
$-2\varpi_1$.

\end{prop}

\setcounter{equation}{0}
\noindent
(B.3)
We now set
$\cE(w)=\cL_\cP(\soc^1_{\ell(w)+1,w})$ with $\cP=G/P$.

\begin{prop}
Let
$x, y\in W^P$.

(i) $\forall i\geq1$,
$\Ext^i_\cP(\cE(x),\cE(y))=0$.

(ii) $\Mod_\cP(\cE(x),\cE(y))\ne0$ iff  $x\geq y$ in the Chevalley-Bruhat order.

\end{prop}

\pf
The assertion is immediate if $s_1s_2s_1\not\in\{x,y\}$. Let us compute the extensions involving $\cE(s_1s_2s_1)$.
We make use of the $P$-structure
on $\soc^1_{4,s_1s_2s_1}$:
\[
\soc^1_{4,s_1s_2s_1}\simeq
(-3\varpi_1)\otimes\nabla(\varpi_1)/
\begin{tabular}{|c|}
\hline
$\nabla^P(-2\varpi_1+\varpi_2)$
\\
\hline
$-\varpi_1$
\\
\hline
\end{tabular}
\simeq
\begin{tabular}{|c|}
\hline
$-2\varpi_1$
\\
\hline
$\nabla^P(-4\varpi_1+\varpi_2)$
\\
\hline
$-3\varpi_1$
\\
\hline
\end{tabular}.
\]

Let $i\in\bbN$.
One has
\begin{align}
\Ext^i_\cP
&
(\cE(e),\cE(s_1s_2s_1))
\simeq
\rH^i(\cP,\cL_\cP(\begin{tabular}{|c|}
\hline
$-2\varpi_1$
\\
\hline
$\nabla^P(-4\varpi_1+\varpi_2)$
\\
\hline
$-3\varpi_1$
\\
\hline
\end{tabular}
)
\\
\notag
&=0
\quad\text{by Bott's theorem
\cite[II.5.5]{J}},
\\
\Ext^i_\cP
&
(\cE(s_1),\cE(s_1s_2s_1))
\simeq
\rH^i(\cP,\cL_\cP(\begin{tabular}{|c|}
\hline
$-\varpi_1$
\\
\hline
$\nabla^P(-3\varpi_1+\varpi_2)$
\\
\hline
$-2\varpi_1$
\\
\hline
\end{tabular}
)
\\
\notag
&=0
\quad\text{likewise},
\\
\Ext^i_\cP
&
(\cE(s_2s_1),\cE(s_1s_2s_1))
\simeq
\rH^i(\cP,\cL_\cP((-\varpi_1)\otimes
\nabla(\varpi_1)/
\begin{tabular}{|c|}
\hline
$\nabla^P(-2\varpi_1+\varpi_2)$
\\
\hline
$-\varpi_1$
\\
\hline
\end{tabular}
))
\\
\notag
&\simeq
\rH^{i+1}(\cP,\cL_\cP(
\begin{tabular}{|c|}
\hline
$\nabla^P(-3\varpi_1+\varpi_2)$
\\
\hline
$-2\varpi_1$
\\
\hline
\end{tabular}
))
=0
\quad\text{likewise},
\\
\Ext^i_\cP
&
(\cE(w^P),\cE(s_1s_2s_1))
\simeq
\rH^i(\cP,\cL_\cP(
\begin{tabular}{|c|}
\hline
$2\varpi_1$
\\
\hline
$\nabla^P(\varpi_2)$
\\
\hline
$\varpi_1$
\\
\hline
\end{tabular}
))
\\
\notag
&\simeq
\delta_{i,0}\{L(2\varpi_1)\oplus
L(\varpi_2)\oplus
L(\varpi_1)\}
\quad\text{by the linkage principle},
\\
\Ext^i_\cP
&
(\cE(s_1s_2s_1), \cE(e))
\simeq
\rH^i(\cP,\cL_\cP(
\begin{tabular}{|c|}
\hline
$3\varpi_1$
\\
\hline
$\nabla^P(\rho)$
\\
\hline
$2\varpi_1$
\\
\hline
\end{tabular}
))
\\
\notag
&\simeq
\delta_{i,0}\{L(3\varpi_1)\oplus
L(\rho)\oplus
L(2\varpi_1)\}
\quad\text{likewise},
\\
\Ext^i_\cP
&
(\cE(s_1s_2s_1), \cE(s_1))
\simeq
\rH^i(\cP,\cL_\cP(
\begin{tabular}{|c|}
\hline
$2\varpi_1$
\\
\hline
$\nabla^P(\varpi_2)$
\\
\hline
$\varpi_1$
\\
\hline
\end{tabular}
))
\\
\notag
&\simeq
\delta_{i,0}\{L(2\varpi_1)\oplus
L(\varpi_2)\oplus
L(\varpi_1)\},
\\
\Ext^i_\cP
&
(\cE(s_1s_2s_1), \cE(s_2s_1))
\simeq
\rH^i(\cP,\cL_\cP(
\begin{tabular}{|c|}
\hline
$\varpi_1$
\\
\hline
$\nabla^P(-\varpi_1+\varpi_2)$
\\
\hline
$\varepsilon$
\\
\hline
\end{tabular}
))
\simeq
\delta_{i,0}\{L(\varpi_1)\oplus
L(e)\},
\end{align}
\begin{align}
\Ext^i_\cP
&
(\cE(s_1s_2s_1), \cE(w^P))
\simeq
\rH^i(\cP,\cL_\cP(
\begin{tabular}{|c|}
\hline
$-\varpi_1$
\\
\hline
$\nabla^P(-3\varpi_1+\varpi_2)$
\\
\hline
$-2\varpi_1$
\\
\hline
\end{tabular}
))
=0,
\\
\Ext^i_\cP
\notag
&
(\cE(s_2s_1s_2s_1),\cE(s_1s_2s_1))
\simeq
\rH^i(\cP,\cL_\cP(\nabla(\varpi_1)/
\begin{tabular}{|c|}
\hline
$\nabla^P(-2\varpi_1+\varpi_2)$
\\
\hline
$-\varpi_1$
\\
\hline
\end{tabular}
)).
\end{align}
There is then a long exact sequence
\begin{multline*}
0\to\nabla(\varpi_1)\to
\Mod_\cP(\cE(s_2s_1s_2s_1),\cE(s_1s_2s_1))
\to
\rH^1(\cP,\cL_\cP(
\nabla^P(-2\varpi_1+\varpi_2)))\to0
\\
\to
\Ext^1_\cP(\cE(s_2s_1s_2s_1),\cE(s_1s_2s_1))
\to
\rH^2(\cP,\cL_\cP(
\nabla^P(-2\varpi_1+\varpi_2)))\to
\dots
\end{multline*}
with
\begin{align*}
\rH^i(\cP,\cL_\cP(
\nabla^P(-2\varpi_1+\varpi_2)))
&\simeq
\rH^i(\cB,\cL(
-2\varpi_1+\varpi_2))
=
\rH^i(\cB,\cL(
s_1\bullet0))
\simeq
\delta_{i,1}L(e),
\end{align*}
and hence by the linkage principle
\begin{align}
\Ext^i_\cP(\cE(s_2s_1s_2s_1),\cE(s_1s_2s_1))
\simeq
\delta_{i,0}\{L(\varpi_1)\oplus
L(e)\}.
\end{align}
 
One has
\begin{align*}
\Ext^i_\cP
&
(\cE(s_1s_2s_1), \cE(s_2s_1s_2s_1))
\simeq
\rH^i(\cP,\cL_\cP(
\begin{tabular}{|c|}
\hline
$\varepsilon$
\\
\hline
$\nabla^P(-2\varpi_1+\varpi_2)$
\\
\hline
$-\varpi_1$
\\
\hline
\end{tabular}
))
\\
&\simeq
\rH^i(\cP,\cL_\cP(
\begin{tabular}{|c|}
\hline
$\varepsilon$
\\
\hline
$\nabla^P(s_1\bullet0)$
\\
\hline
\end{tabular}
)),
\end{align*}
and hence a long exact sequence
\begin{multline*}
0\to
\Mod_\cP
(\cE(s_1s_2s_1), \cE(s_2s_1s_2s_1))
\to
L(e)\to
\rH^1(\cB,\cL(
s_1\bullet0))\to
\\
\Ext^1_\cP
(\cE(s_1s_2s_1), \cE(s_2s_1s_2s_1))
\to0\to
\dots
\end{multline*}
with
$\rH^i(\cB,\cL(
s_1\bullet0))\simeq\delta_{i,1}L(e)$.
On the other hand,
\begin{align*}
\Mod_\cP
&
(\cE(s_1s_2s_1), \cE(s_2s_1s_2s_1))
\simeq
\Mod_\cP
(\cL_\cP(\nabla(\varpi_1)/
\begin{tabular}{|c|}
\hline
$\nabla^P(-2\varpi_1+\varpi_2)$
\\
\hline
$-\varpi_1$
\\
\hline
\end{tabular}
), \cO_\cP)
\\
&\leq
\Mod_\cP
(\cL_\cP(\nabla(\varpi_1)
), \cO_\cP)
\simeq
L(\varpi_1).
\end{align*}
It follows that
\begin{align}
\Ext^\bullet_\cP
(\cE(s_1s_2s_1), \cE(s_2s_1s_2s_1))=0.
\end{align}

Finally,
\begin{align}
\Ext^i_\cP
&
(\cE(s_1s_2s_1), \cE(s_1s_2s_1))
\\
\notag&
\simeq
\Ext^i_\cP
(\cE(s_1s_2s_1), \cL_\cP((-3\varpi_1)\otimes
\nabla(\varpi_1)/
\begin{tabular}{|c|}
\hline
$\nabla^P(-2\varpi_1+\varpi_2)$
\\
\hline
$-\varpi_1$
\\
\hline
\end{tabular}
))
\\
\notag
&\simeq
\Ext^{i+1}_\cP
(\cE(s_1s_2s_1), \cL_\cP((-3\varpi_1)\otimes
\begin{tabular}{|c|}
\hline
$\nabla^P(-2\varpi_1+\varpi_2)$
\\
\hline
$-\varpi_1$
\\
\hline
\end{tabular}
))
\quad\text{as}
\\
\notag
&\hspace*{3cm}\text{
$\Ext^{\bullet}_\cP
(\cE(s_1s_2s_1), \cL_\cP((-3\varpi_1)\otimes
\nabla(\varpi_1)))\simeq$}
\\
\notag
&\hspace*{3cm}\text{
$\Ext^{\bullet}_\cP
(\cE(s_1s_2s_1), \cE(s_2s_1s_2s_1))
\otimes
\nabla(\varpi_1)=0$
by
(10)}
\\
\notag
&\simeq
\Ext^{i+1}_\cP
(\cE(s_1s_2s_1), \cL_\cP(
\begin{tabular}{|c|}
\hline
$\nabla^P(-5\varpi_1+\varpi_2)$
\\
\hline
$-4\varpi_1$
\\
\hline
\end{tabular}
))
\\
\notag
&\simeq
\Ext^{i+1}_\cP
(\cE(s_1s_2s_1), \cL_\cP(
\nabla^P(-5\varpi_1+\varpi_2)))
\quad\text{by (8)}
\\
\notag
&\simeq
\Ext^{i+1}_\cP
(\cE(\begin{tabular}{|c|}
\hline
$\varpi_1$
\\
\hline
$\nabla^P(-\varpi_1+\varpi_2)$
\\
\hline
$\varepsilon$
\\
\hline
\end{tabular}
), \cL_\cP(
\nabla^P(-2\varpi_1+\varpi_2)))
\\
\notag
&\simeq
\rH^{i+1}(\cP,
\cL_\cP(\begin{tabular}{|c|}
\hline
$\varepsilon$
\\
\hline
$\nabla^P(-2\varpi_1+\varpi_2)$
\\
\hline
$-\varpi_1$
\\
\hline
\end{tabular}
\otimes
\nabla^P(-2\varpi_1+\varpi_2)))
\\
\notag
&\simeq
\rH^{i+1}(\cP,
\cL_\cP(\begin{tabular}{|c|}
\hline
$\nabla^P(-2\varpi_1+\varpi_2)$
\\
\hline
$\nabla^P(-4\varpi_1+2\varpi_2)$
\\
\hline
$\nabla^P(-\varpi_1)$
\\
\hline
$\nabla^P(-3\varpi_1+\varpi_2)$
\\
\hline
\end{tabular}
))
\simeq
\rH^{i+1}(\cP,
\cL_\cP(\nabla^P(-2\varpi_1+\varpi_2)))
\\
\notag
&=
\rH^{i+1}(\cP,
\cL_\cP(\nabla^P(s_1\bullet0)))
\simeq
\delta_{i+1,1}L(e)
\simeq
\delta_{i,0}L(e).
\end{align}

\setcounter{equation}{0}
\noindent
(B.4)
We show that
$\cE(w)$, $w\in W^P$, Karoubian generate
$\rD^b(\coh\cP)$ as in (A.9).
Let
$\hat\cE=\langle\cE(w)|w\in W^P\rangle$
denote the triangulated subcategory of
$\rD^b(\coh\cP)$ Karoubian generated by the
$\cE(w)$, $w\in W^P$.

As
$\Lambda_P=\bbZ\varpi_1$, it is enough to show that
all $\cL_\cP(n\varpi_1)\in\hat\cE$,
$n\in\bbZ$.
For that
we may transfer to
$\cB=G/B$ and show that
all $\cL(n\varpi_1)\in\pi^*\hat\cE$,
$n\in\bbZ$.
Put
$\tilde\cE=\pi^*\hat\cE$.
For our purpose
we may also assume that,
whenever $\cL(\nabla^P(M))\in\tilde\cE$ for a $P$-module
$M$, $\cL(M)\in\tilde\cE$, and vice versa.
In particular,
if
$
\cL(\nabla^{P}(n\varpi_1+\varpi_2))
\in\tilde\cE$, $n\in\bbZ$, then
$\cL(n\varpi_1+\varpi_2)
\in\tilde\cE$,
and hence also
$\cL((n+3)\varpi_1-\varpi_2)
\in\tilde\cE$.

Now, using $\cE(s_1s_2s_1),\cE(s_2s_1)$ and $\cE(s_2s_1s_2s_1)$, we see that
$\cL_\cP(\nabla^P(-4\varpi_1+\varpi_2))\in\hat\cE$, and hence
$\cL(-4\varpi_1+\varpi_2)$ and $\cL(-\rho)\in\tilde\cE$. 
As
\begin{align*}
\hat\cE
&\ni\nabla(\varpi_1)\otimes
\cL_\cP(-2\varpi_1)
\simeq
\cL_\cP(\begin{tabular}{|c|}
\hline
$\varpi_1$
\\
\hline
$\nabla^{P}(-\varpi_1+\varpi_2)$
\\
\hline
$\varepsilon$
\\
\hline
$\nabla^{P}(-2\varpi_1+\varpi_2)$
\\
\hline
$-\varpi_1$
\\
\hline
\end{tabular})
\otimes
(-2\varpi_1))
\simeq
\cL_\cP(\begin{tabular}{|c|}
\hline
$-\varpi_1$
\\
\hline
$\nabla^{P}(-3\varpi_1+\varpi_2)$
\\
\hline
$-2\varpi_1$
\\
\hline
$\nabla^{P}(-4\varpi_1+\varpi_2)$
\\
\hline
$-3\varpi_1$
\\
\hline
\end{tabular}),
\end{align*}
$\cL_\cP(\nabla^{P}(-3\varpi_1+\varpi_2))
\in\hat\cE$,
and hence
$\cL(-3\varpi_1+\varpi_2))$
and
$\cL(-\varpi_2)\in\tilde\cE$.
As
\[
\hat\cE\ni\nabla(\varpi_1)\otimes
\cL_\cP(-\varpi_1)
\simeq
\cL_\cP(\begin{tabular}{|c|}\hline$\varepsilon$\\\hline$\nabla^{P}(-2\varpi_1+\varpi_2)$\\\hline$-\varpi_1$\\\hline$\nabla^{P}(-3\varpi_1+\varpi_2)$\\\hline$-2\varpi_1$\\\hline
\end{tabular}),
\]
$\cL(-2\varpi_1+\varpi_2)$
and
$\cL(\varpi_1-\varpi_2)
\in\tilde\cE$.
As
\[
\hat\cE\ni\nabla(\varpi_1)\otimes
\cL_\cP(-3\varpi_1)
\simeq
\cL_\cP(\begin{tabular}{|c|}
\hline
$-2\varpi_1$
\\
\hline
$\nabla^{P}(-4\varpi_1+\varpi_2)$
\\
\hline
$-3\varpi_1$
\\
\hline
$\nabla^{P}(-5\varpi_1+\varpi_2)$
\\
\hline
$-4\varpi_1$
\\
\hline
\end{tabular}),
\]
$\cL(-5\varpi_1+\varpi_2))$
and
$\cL(-2\varpi_1-\varpi_2)
\in\tilde\cE$.
As
\begin{align*}
\hat\cE
&\ni\nabla(\varpi_2)\otimes
\cL_\cP(-2\varpi_1)
\\
&\simeq
\cL_\cP(\begin{tabular}{|c|}
\hline
$\nabla^{P}(\varpi_2)$
\\
\hline
$\varpi_1$
\\
\hline
$\nabla^{P}(-\varpi_1+\varpi_2)$
\\
\hline
$\nabla^{P}(-3\varpi_1+2\varpi_2)\oplus
\varepsilon$
\\
\hline
$\nabla^{P}(-2\varpi_1+\varpi_2)$
\\
\hline
$-\varpi_1$
\\
\hline
$\nabla^{P}(-3\varpi_1+\varpi_2)$
\\
\hline
\end{tabular})
\otimes
(-2\varpi_1))
\simeq
\cL_\cP(\begin{tabular}{|c|}
\hline
$\nabla^{P}(-2\varpi_1+\varpi_2)$
\\
\hline
$-\varpi_1$
\\
\hline
$\nabla^{P}(-3\varpi_1+\varpi_2)$
\\
\hline
$\nabla^{P}(-5\varpi_1+2\varpi_2)
\oplus
(-2\varpi_1)$
\\
\hline
$\nabla^{P}(-4\varpi_1+\varpi_2)$
\\
\hline
$-3\varpi_1$
\\
\hline
$\nabla^{P}(-5\varpi_1+\varpi_2)$
\\
\hline
\end{tabular}),
\end{align*}
$
\cL(-5\varpi_1+2\varpi_2))$
and $\cL(\varpi_1-2\varpi_2)
\in\tilde\cE$.
As
$\nabla(\varpi_1)\otimes
\cL(-\varpi_2)\in\tilde\cE$,
$\cL(2\varpi_1-2\varpi_2)\in\tilde\cE$.
As
$\nabla(\varpi_1)\otimes
\cL(-\rho)\in\tilde\cE$,
$\cL(-2\varpi_2)\in\tilde\cE$.
As
$\nabla(\varpi_1)\otimes
\cL(-3\varpi_1+\varpi_2)\in\tilde\cE$,
$\cL(-4\varpi_1+2\varpi_2)\in\tilde\cE$.
As
$\nabla(\varpi_2)\otimes
\cL(-\varpi_1)\in\tilde\cE$,
$\cL(\nabla^{P}(-\varpi_1+\varpi_2))\in\tilde\cE$.
Then
$\cL(-\varpi_1+\varpi_2)\in\tilde\cE$, and hence also
$\cL(2\varpi_1-\varpi_2)\in\tilde\cE$.
As $\nabla(\varpi_1)\otimes\cO_\cB\in\tilde\cE$, $\cL(\varpi_1)\in\tilde\cE$.
As
$\nabla(\varpi_1)\otimes
\cL(-4\varpi_1+\varpi_2)\in\tilde\cE$,
$\cL(-6\varpi_1+2\varpi_2)\in\tilde\cE$.
As
$\nabla(\varpi_2)\otimes
\cL(-3\varpi_1)\in\tilde\cE$,
$\cL(\nabla^{P}(-6\varpi_1+\varpi_2))\in\tilde\cE$.
Then
$\cL(-6\varpi_1+\varpi_2)\in\tilde\cE$, and hence also
$\cL(-3\varpi_1-\varpi_2)\in\tilde\cE$.
As
$\nabla(\varpi_1)\otimes\cL(-4\varpi_1)\in\tilde\cE$,
$\cL(-5\varpi_1)\in\tilde\cE$.
Thus
\[
\cL(k\varpi_1)\in\tilde\cE
\quad\forall k\in[-5,1].
\]
As
$\dim\nabla(\varpi_1)=7$,
one now obtains all
$\cL(n\varpi_1)\in\tilde\cE$,
$n\in\bbZ$,
from the exact sequence
\[
0\to
\cL(-7\varpi_1)
\to
\cL(-6\varpi_1)^{\oplus\binom{7}{6}}
\to
\cL(-5\varpi_1)^{\oplus\binom{7}{5}}
\to
\dots
\to
\cL(-\varpi_1)^{\oplus\binom{7}{1}}
\to
\cO_\cB
\to
0.
\]

Thus, we have obtained
\begin{thm}
Assume $p\geq11$.
The decomposition of $F_*\cO_\cP$ into indecomposables is given by
\begin{multline*}
\cE(e)
\oplus
\cE(s_1)\otimes L(s_1)
\oplus
\cE(s_2s_1)\otimes\{L(s_2s_1)\oplus L(e)\oplus
L(w^P)\}
\oplus
\cE(s_1s_2s_1)\otimes L(s_1s_2s_1)
\\
\oplus
\cE(s_2s_1s_2s_1)\otimes\{L(s_2s_1s_2s_1)\oplus
L(s_1)\}
\oplus
\cE(w^P)\otimes\{L(w^P)\oplus L(e)\}
\end{multline*}
with the $\cE(w)$, $w\in W^P$, 
forming a Karoubian complete strongly exceptional poset 
such that
$\Mod_\cP(\cE(x),\cE(y))\ne0$
iff $x\geq y$ in the Chevalley-Bruhat order.
\end{thm}

\noindent
{\bf C.}
We also append an explicit imbedding 
of $G$ into $\rSO_7(\Bbbk)$,
which is essentially the same as
H\'{e}e's \cite[13.6]{He}.
We will allow $\Bbbk$ to be any algebraically closed field of characteristic not 2.
The author is grateful to Tanisaki and Testerman for references.

\setcounter{equation}{0}
\noindent
(C.1)
We show first the imbedding of
a $\bbZ$-form 
$\frg_\bbZ$
of 
the Lie algebra 
of $G$ into 
a $\bbZ$-form
$\frg'_\bbZ$
of the Lie algebra  of $\rSO_7(\Bbbk)$.
For that we recall imbeddings of $\bbQ$-Lie algebras,
a classical result called the principle of triality,
using folding
\cite{T}.
Let
$\tilde\frg$ be the Lie algebra of $\rSO_8(\bbQ)$
with Dynkin diagram
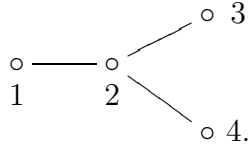
\[
\xymatrix@R0ex{
\ & \ & 
\circ
& \hspace{-1.5cm} 3
\\
\\
\circ
\ar@{-}[r]
&
\circ
\ar@{-}[uur]\ar@{-}[ddr]
\\
1
&
2
& &
\\
\ & \ & 
\circ & \hspace{-1.5cm}4. 
}
\]
Let $A$ be the associated Cartan matrix, and let 
$\tilde e_i,\tilde h_i, \tilde f_i$, $i\in[1,4]$,
be the 
standard Chevalley generators of $\tilde\frg$
such that
$[\tilde e_i,\tilde f_j]=\delta_{i,j}\tilde h_i$,
$[\tilde h_i, \tilde e_j]=
A_{ij}\tilde e_i$,
$[\tilde h_i, \tilde f_j]=-
A_{ij}\tilde f_j$
$\forall i, j\in[1,4]$.
Let
$\sigma$ ba an automorphism of the Dynkin diagram;
$
A_{\sigma(i)\sigma(j)}=
A_{ij}$ $\forall i,j\in[1,4]$.
By the same letter $\sigma$ we let it also denote the induced automorphism of $\tilde\frg$ such that
$\tilde e_i\mapsto\tilde e_{\sigma(i)}$,
$\tilde f_i\mapsto\tilde f_{\sigma(i)}$,
$\forall i\in[1,4]$.
Let
$\tilde\frg^\sigma=\{x\in\tilde\frg|\sigma(x)=x\}$ be the fixed point subalgebra of $\tilde\frg$
under $\sigma$.
Let
$\cO$ be a $\langle\sigma\rangle$-orbit in the index set $[1,4]$.
We divide into the following 2 cases.

Case 1:
Either
$|\cO|=1$ or $A_{ji}=0$ for any distinct $i,j\in\cO$,

Case 2: $\cO=\{i,j\}$ with $i\ne j$ such that
$A_{ij}=-1=A_{ji}$.

\noindent
For each orbit
$\cO$ define elements of $\tilde\frg^\sigma$ by
\[
\tilde h_\cO=
\begin{cases}
\sum_{i\in\cO}\tilde h_i
&\text{in Case 1},
\\
2\sum_{i\in\cO}\tilde h_i
&\text{in Case 2},
\end{cases}
\quad
\tilde e_\cO=
\sum_{i\in\cO}\tilde e_i,
\quad
\tilde f_\cO=
\begin{cases}
\sum_{i\in\cO}\tilde f_i
&\text{in Case 1},
\\
2\sum_{i\in\cO}\tilde f_i
&\text{in Case 2}.
\end{cases}
\]

Now, let
$\frg_\bbQ=\frg_\bbZ\otimes_\bbZ\bbQ$ 
and let
$e_1, e_2, f_1,f_2$
be the Chevalley generators
corresponding to the simple
roots
$\alpha_1$ and $\alpha_2$.
Taking $\sigma$ of order 3, one obtains from \cite[Th. B.4]{T}
an isomorphism of 
$\bbQ$-Lie algebras
$\theta_1:\frg_\bbQ\to\tilde\frg^\sigma$ such that
\[
e_1\mapsto\tilde e_1+\tilde e_3+\tilde e_4,
\quad
e_2\mapsto\tilde e_2,
\quad
f_1\mapsto\tilde f_1+\tilde f_3+\tilde f_4,
\quad
f_2\mapsto\tilde f_2.
\]

Let next
$V$ be a 7-dimensional $\Bbbk$-linear space with
basis $v_1, v_2, v_3, v_0, v_{-3}, v_{-2},v_{-1}$ equipped with a quadratic form
$Q(\sum_{k=-3}^3\xi_kv_k)=
\xi_0^2+
\sum_{k=1}^3\xi_k\xi_{-k}$ 
for $\xi_i,\xi_{-i}\in\Bbbk$, $i\in[0,3]$.
Thus the associated Gram matrix is
\[
[\!(\bbB(v_i,v_j)=q(v_i+v_j)-q(v_i)-q(v_j))\!]
=
\left(\begin{array}{ccc|c|ccc}
&&&&&&1
\\
&&&&&1&
\\
&&&&1&&
\\
\hline
&&&2&&&
\\
\hline
&&1&&&&
\\
&1&&&&&
\\
1&&&&&&
\end{array}
\right).
\]
We regard $G'=\rSO(V;\bbB)=\{g\in\SL(V)|\bbB(gv,gv')=\bbB(v,v')\ \forall
v,v'\in V\}$
as our orthogonal group
$\rSO_7(\Bbbk)$.
Let
$T'=
\{\diag(\zeta_1,\zeta_2,\zeta_3,1,\zeta_3^{-1},\zeta_2^{-1},\zeta_1^{-1})|\zeta_1,\zeta-2\in\Bbbk^\times\}$ be a maximal torus of $G'$ with simple coroots
$\alpha_1^\vee=
\varepsilon_1^\vee-\varepsilon_2^\vee,
\alpha_2^\vee=
\varepsilon_2^\vee-\varepsilon_{3}^\vee$,
and $\alpha_{3}^\vee=2\varepsilon_{3}^\vee$,
where
$\varepsilon_1^\vee : \zeta\mapsto
\diag(\zeta, 1,1,
1,1,1,\zeta^{-1})$,
$\varepsilon_2^\vee : \zeta\mapsto
\diag(1,\zeta, 1,1,1,1,\zeta^{-1},1)
$,
and
$\varepsilon_3^\vee : \zeta\mapsto
\diag(1,1,\zeta, 1,\zeta^{-1},1,1)
$.
If
$\varepsilon_k:$
$\diag(\zeta_1,\zeta_2,\zeta_3,\zeta_0, \zeta_3^{-1},\zeta_2^{-1},\zeta_1^{-1})\mapsto\zeta_k$,
$k\in[1,3]$,
the corresponding simple roots are 
$\alpha'_1=\varepsilon_1-\varepsilon_2,
\alpha'_2=\varepsilon_2-\varepsilon_3$,
and $\alpha'_3=\varepsilon_3$.
If we let
$E$ denote the identity matrix and $E(i,j)$,
$i,j\in[-3,3]$,
denote the square matrix of degree 7
with 1 at the $(i,j)$-th entry and 0 elsewhere,
the root subgroups of $G'$
are
given by,
for
$i,j\in[1,3]$
with
$i<j$,
\begin{align*}
U_{\varepsilon_i-\varepsilon_j}
&
=
\{
E+\xi(E(i,j)-E(-j,-i))
\mid
\xi\in\Bbbk\},
\\ 
U_{-\varepsilon_i+\varepsilon_j}
&
=
\{
E+\xi(E(j,i)-E(-i,-j))
\mid
\xi\in\Bbbk\},
\\
U_{\varepsilon_i+\varepsilon_j}
&
=
\{
E+\xi(E(i,-j)-E(j,-i))
\mid
\xi\in\Bbbk\},
\\
U_{-\varepsilon_i-\varepsilon_j}
&
=
\{
E+a(E(-j,i)-E(-i,j))
\mid
a\in\Bbbk\},
\end{align*}
and, for $k\in[1,3]$,
\begin{align*}
U_{\varepsilon_k}
&=
\{
E+\xi(2E(k,0)-E(0,-k))-\xi^2E(k,-k)
\mid
\xi\in\Bbbk\},
\\
U_{-\varepsilon_k}
&=
\{
E+\xi(E(0,k)-2E(-k,0))-\xi^2E(-k,k)
\mid
\xi\in\Bbbk\}.
\end{align*}
\[
\xymatrix@R0ex{
\varepsilon_1-\varepsilon_2
&
\varepsilon_2-\varepsilon_3
&
\varepsilon_{3}
\\
\ar@{-}[r]
&
\circ
\ar@{=>}[r]
&
\circ
\\
1
&
2
&
3
}
\]
Then
$e'_1=E(1,2)-E(-2,-1),
e'_2=E(2,3)-E(-3,-2),
e'_3=2E(3,0)-E(0,-3)$,
and $f'_1=E(2,1)-E(-1,-2),
f'_2=E(3,2)-E(-2,-3),
f'_3=E(0,3)-2E(-3,0)$
form Chevalley generators of the $\bbQ$-Lie algebra
$\frg'_\bbQ=\frg'_\bbZ\otimes_\bbZ\bbQ$.
If $\sigma$ is of order 2, one obtains by \cite[Th. B.4]{T}
an isomorphism of Lie algebras
$\theta_2:\frg'_\bbQ\to\tilde\frg^\sigma$ such that
\[
e'_1\mapsto\tilde e_1,
\quad
e'_2\mapsto\tilde e_2,
\quad
e'_3\mapsto\tilde e_3+\tilde e_4,
\quad
f'_1\mapsto\tilde f_1,
\quad
f'_2\mapsto\tilde f_2,
\quad
f'_3\mapsto\tilde f_3+\tilde f_4.
\]
It follows that
$\theta_1$ factors through 
$\theta_2$
to yield an imbedding
$\theta_\bbQ:\frg_\bbQ\hookrightarrow\frg'_\bbQ$ of $\bbQ$-Lie algebras.

In $\frg_\bbQ$
one can take
along with $e_i,f_i$ and $[e_i,f_i]$, $i\in[1,2]$,
\begin{align}
e_{\alpha_1+\alpha_2}
&=
[e_1,e_2],
\ 
e_{2\alpha_1+\alpha_2}=\frac{1}{2}[e_1,e_{\alpha_1+\alpha_2}],
\ 
e_{3\alpha_1+\alpha_2}
=
\frac{1}{3}
[e_1,e_{2\alpha_1+\alpha_2}],
\\ 
e_{3\alpha_1+2\alpha_2}
\notag&=
[e_2,e_{3\alpha_1+\alpha_2}],
\\
f_{\alpha_1+\alpha_2}
\notag&=
-[f_1,f_2],
\ 
f_{2\alpha_1+\alpha_2}
=
-\frac{1}{2}
[f_1,f_{\alpha_1+\alpha_2}],
\ 
f_{3\alpha_1+\alpha_2}
=
-\frac{1}{3}[f_1,f_{2\alpha_1+\alpha_2}],
\\ 
f_{3\alpha_1+2\alpha_2}
\notag&=
-[f_2,f_{3\alpha_1+\alpha_2}]
\end{align}
to form a Chevalley basis of $\frg_\bbZ$
\cite[Th. B.2.1]{Ca}.
Also,
along with 
$e'_i,f'_i$ and $[e'_i,f'_i]$, $i\in[1,3]$, 
\begin{multline}
\{E(i,j)-E(-j,-i),
E(j,i)-E(-i,-j),
E(i,-j)-E(j,-i),
E(-j,i)-E(-i,j)
\\
|1\leq i<j\leq3\}
\sqcup
\{
2E(k,0)-E(0,-k),
E(0,k)-2E(-k,0)
|
k\in[1,3]\}
\end{multline}
form a Chevalley basis of
$\frg'_\bbZ$
\cite[11.2.4]{Ca}.
Then
\begin{align}
\theta(e_{\alpha_1+\alpha_2})
&=
(E(1,3)-E(-3,-1))-
(2E(2,0)-E(0,-2)),
\\
\theta(e_{2\alpha_1+\alpha_2})
\notag&=
-(2E(1,0)-E(0,-1))-
(E(2,-3)-E(3,-2)),
\\
\theta(e_{3\alpha_1+\alpha_2})
\notag&=
-(E(1,-3)-E(3,-1)),
\\
\theta(e_{3\alpha_1+2\alpha_2})
\notag&=
-(E(1,-2)-E(2,-1)),
\\
\theta(f_{\alpha_1+\alpha_2})
\notag&=
-(E(0,2)-2E(-2,0))+
(E(3,1)-E(-1,-3)),
\\
\theta(f_{2\alpha_1+\alpha_2})
\notag&=
-(E(0,1)-2E(-1,0))-
(E(-3,2)-E(-2,3)),
\\
\theta(f_{3\alpha_1+\alpha_2})
\notag&=
-(E(-3,1)-E(-1,3)),
\\
\theta(f_{3\alpha_1+2\alpha_2})
\notag&=
-(E(-2,1)-E(-1,2)).
\end{align}
Thus,

\begin{prop}
There is an imbedding of Lie algebras
$\theta_\bbZ: \frg_\bbZ\to\frg'_\bbZ$ such that
\[
e_1\mapsto e'_1+e'_3,
\quad
e_2\mapsto e'_2,
\quad
f_1\mapsto
f'_1+f'_3,
\quad
f_2\mapsto
f'_2.
\]

\end{prop}


\setcounter{equation}{0}
\noindent
(C.2)
Using the representation 
$\theta_\bbZ$
of $\frg_\bbZ$ on
$\bbZ^{\oplus_7}$,
we exponentiate to obtain a realization of $G$ in $\GL_7(\Bbbk)$
factoring through $G'$
\cite{St}, which is
essentially the same as
\cite[13..6]{He}.
Let $y_\alpha$, $\alpha\in R$, denote the root vectors of the Chevalley basis (C.1) of $\frg_\bbZ$, and put
$y'_\alpha=\theta(y_\alpha)$,
$x_\alpha(\xi)=
\exp(\xi y'_\alpha)$.
\begin{prop}
One has an imbedding of algebraic
groups
$G\to\rSO_7(\Bbbk)$ with
the root subgroups
given by
$U_\alpha=
\{
x_\alpha(\xi)|
\xi\in\Bbbk\}$,
$\alpha\in R$.

\end{prop}

\setcounter{equation}{0}
\noindent
(C.3)
Explicitly,
the root subgroups 
$U_\alpha$, $\alpha\in R$, of $G$ are realized in $G'$ as follows:
\begin{align}
x_{\alpha_1}(\xi)
&=
E+\xi(E(1,2)+2E(3,0)-E(0,-3)-E(-2,-1))-\xi^2E(3,-3),
\\
x_{\alpha_2}(\xi)
\notag&=
E+\xi(E(2,3)-E(-3,-2)),
\\
x_{\alpha_1+\alpha_2}(\xi)
\notag&=
E+
\xi(E(1,3)-2E(2,0)+E(0,-2)-E(-3,-1))
-\xi^2E(2,-2),
\\
x_{2\alpha_1+\alpha_2}(\xi)
\notag&=
E+
\xi(-2E(1,0)-E(2,-3)+E(3,-2)+E(0,-1))
-\xi^2E(1,-1),
\\
x_{3\alpha_1+\alpha_2}(\xi)
\notag&=
E+\xi(-E(1,-3)+E(3,-1)),
\\
x_{3\alpha_1+2\alpha_2}(\xi)
\notag&=
E+\xi(-E(1,-2)+E(2,-1)),
\\
x_{-\alpha_1}(\xi)
\notag&=
E+\xi(E(2,1)+E(0,3)-2E(-3,0)-E(-1,-2))
-\xi^2E(-3,3),
\\
x_{-\alpha_2}(\xi)
\notag&=
E+\xi(E(3,2)-E(-2,-3)),
\\
x_{-\alpha_1-\alpha_2}(\xi)
\notag&=
E+\xi(E(3,1)-E(0,2)+2E(-2,0)-E(-1,-3))
-\xi^2E(-2,2),
\\
x_{-2\alpha_1-\alpha_2}(\xi)
\notag&=
E+\xi(-E(0,1)-E(-3,2)+E(-2,3)+2E(-1,0))
-\xi^2E(-1,1),
\\
x_{-3\alpha_1-\alpha_2}(\xi)
\notag&=
E+\xi(-E(-3,1)+E(-1,3)),
\\
x_{-3\alpha_1-2\alpha_2}(\xi)
\notag&=
E+\xi(-E(-2,1)+E(-1,2)).
\end{align}

In particular,
if we let
$x_{\alpha'_i}(\xi)=\exp(\xi
e'_i)$ and
$x_{-\alpha'_i}(\xi)=\exp(\xi
f'_i)$,
$i\in[1,3]$,
\begin{alignat}{2}
x_{\alpha_1}(\xi)
&=
x_{\alpha'_1}(\xi)x_{\alpha'_3}(\xi),\ 
&
x_{\alpha_2}(\xi)
&=
x_{\alpha'_2}(\xi),
\\ 
x_{-\alpha_1}(\xi)
\notag&=
x_{-\alpha'_1}(\xi)x_{-\alpha'_3}(\xi),\quad 
&
x_{-\alpha_2}(\xi)
&=
x_{-\alpha'_2}(\xi).
\end{alignat}
$\forall\alpha\in R^s$
$\forall\zeta\in\Bbbk^\times$, set
after 
\cite[p. 43]{St}
\begin{align}
n_\alpha(\zeta)=x_\alpha(\zeta)
x_{-\alpha}(-\zeta^{-1})x_\alpha(\zeta),
\quad
\alpha^\vee(\zeta)=
n_\alpha(\zeta)n_\alpha(-1).
\end{align}
As $[e'_1,e'_3]=0$ in $\frg'$, one has in $G'$
\begin{align}
\alpha_1^\vee(\zeta)
&=
{\alpha'_1}^\vee(\zeta){\alpha'_3}^\vee(\zeta)
=
\diag(\zeta\ \zeta^{-1}\ \zeta^2\ 1\ \zeta^{-2}
\ \zeta\ \zeta^{-1}),
\\
\alpha_2^\vee(\zeta)
\notag&=
{\alpha'_2}^\vee(\zeta)
=
\diag(1\ \zeta\ \zeta^{-1}\  1\ \zeta
\ \zeta^{-1}\ 1).
\end{align}
It follows that 
the fundamental weights
$\varpi'_1=\varepsilon_1$,
$\varpi'_2=\varepsilon_1+\varepsilon_2$, and 
$\varpi'_3=\frac{1}{2}(\varepsilon_1+\varepsilon_2+\varepsilon_3)$ for $T'$
read as $T$-characters
\begin{align}
\varpi'_1|_T=\varpi_1=\varpi'_3|_T,
\quad
\varpi'_2|_T=\varpi_2,
\end{align}
and that
the $T$-weight $\wt(v_k)$
of $v_k$, $k\in[-3,3]$, are 
\begin{alignat}{2}
\wt(v_1)
&=
\varpi_1=-\wt(v_{-1}),
&
\wt(v_2)
&=
-\varpi_1+\varpi_2=
-\wt(v_{-2}),
\\
\wt(v_3)
\notag&=
2\varpi_1-\varpi_2
=
-\wt(v_{-3}),
\quad
&\wt(v_0)
&=
0.
\end{alignat}

\setcounter{equation}{0}
\noindent
(C.4)
{\bf Remarks:}
(i)
The realization of $G$ in $\GL_7(\Bbbk)$ as above holds, of course, over any field $\Bbbk$.

(ii)
The ambient space $V$ as a $G$-module affords
$\nabla(\varpi_1)$, which remains simple over any field of odd characteristic.
In characteristic 2, however,
$G$ stabilizes $\Bbbk v_0$, and hence
$V$ is rather isomorphic to the Weyl module
$\Delta(\varpi_1)$
over any field.
By a base change 
and modulo sign changes in the Chevalley basis
the presentation of the root subgroups coincides
with the one given in \cite[p. 43]{Te}.

\setcounter{equation}{0}
\noindent
(C.5)
Let
$B'$ be the Borel subgroup of $G'$ consisting of lower triangular matrices. 
By the unicity of parabolic subgroups
the stabilizer in $G'$ (resp. $G$)
of the line
$\Bbbk v_{-1}$ is the standard
parabolic subgroup
$P'_{\{\alpha'_2,\alpha'_3\}}$ 
(resp. $P_{\alpha_2}$)
of $G'$
(resp. $G$).
One thus obtains an injective morphism
$\phi:G/P_{\alpha_2}\to G'/P'_{\{\alpha_2,\alpha_3\}}$.
If $B^+$
(resp. ${B'}^+$) is the Borel subgroup of $G$
(resp. $G'$) opposite to
$B$
(resp. $B'$), 
using $\rd({\eta'}^{-1}):\frg \hookrightarrow\frg'$,
one sees that
$\phi$ induces an isomorphism
$B^+P_{\alpha_2}/P_{\alpha_2}\to{B'}^+P'_{\{\alpha_2,\alpha_3\}}/P'_{\{\alpha_2,\alpha_3\}}$.
As the latter is open in $G'/P'_{\{\alpha_2,\alpha_3\}}$
and as
$G/P_{\alpha_2}$ is complete,
$\phi$ must itself be an isomorphism.
Likewise
the stabilizer of $\Bbbk
v_{-2}\oplus\Bbbk v_{-1}$ in $G$ is
$P=P_{\alpha_1}$.
Thus

\begin{cor}
There is an isomorphism of varieties $G/P_{\alpha_2}\simeq
G'/P'_{\{\alpha_2,\alpha_3\}}$ and a closed imbedding
$G/P\hookrightarrow\Gr(2,7)$.

\end{cor}


\begin{thebibliography}{AAAAA}

\bibitem[AbK]{AbK}  Abe, N. and Kaneda, M., {\it On the structure of parabolically induced $G_1T$-Verma modules}, JIM Jussieu {\bf 14} Issue 01 (2015), 185-220

\bibitem[A86]{A86} Andersen, H.H., {\it An inversion formula for the Kazhdan-Lusztig polynomials for affine Weyl groups}, Adv. Math. {\bf 60}
(1986), 125-153






\bibitem[AK89]{AK89}  Andersen, H.H. and Kaneda M., {\it Loewy series of modules for the first Frobenius kernel in a reductive algebraic group}, Proc. LMS (3) {\bf 59} (1989), 74--98

\bibitem[AK00]{AK00}  Andersen, H.H. and Kaneda M., {\it On the $D$-affinity of the flag variety in type $B\sb 2$}, Manuscripta Math. {\bf 103} (2000), no. 3, 393--399














\bibitem[Ber]{B}  Berthelot, P., {\it $\cD$-modules arithm\'etiques I. Op\'erateurs diff\'erentiels de niveau fini}, Ann. scient. \'Ec. Norm. Sup. {\bf 29} (1996), 185-272.

\bibitem[BMR02]{BMR02}  Bezrukavnikov, R., Mirkovic, I. and Rumynin, D., {\it Localization of modules for a semisimple Lie algebra in prime characteristic}, arXiv:math.RT/0205144v1





\bibitem[BMR]{BMR}  Bezrukavnikov, R., Mirkovic, I. and Rumynin, D., {\it Singular localization and intertwining functors  for reductive Lie algebras in prime characteristic}, Nagoya Math. J. {\bf 184} (2006), 1--55


\bibitem[Ca]{Ca} Carter, R., Simple Groups of Lie Type, London etc. 1972 (Wiley)





\bibitem[DG]{DG}
Donkin, S. and Geranios, H.,
{\it 
The cohomology of line bundles for the symplectic group of degree 4 in characteristic 2},
unpublished manuscript



\bibitem[GK]{GK} Gros, M. and Kaneda, M., {\it Contraction par Frobenius
et modules de Steinberg}, preprint


\bibitem[EGA]{EGA} Grothendieck, A. and Dieudonn\'e, J., \lq\lq \'El\'ements de G\'eom\'etrie Alg\'ebrique IV", Pub. Math. no. 24, IHES 1965.


\bibitem[Haa]{Haa} Haastert, B., {\it \"{U}ber Differentialoperatoren und $\Bbb D$-Moduln in positiver Charakteristik}, Manusc. Math. 58 (1987), 385--415








\bibitem[HKR]{HKR} Hashimoto Y., Kaneda M. and Rumynin, D., {\it On localization of $\bar D$-modules}, in ``Representations of Algebraic Groups, Quantum Groups, and Lie Algebras," Contemp. Math. 413 (2006), 43-62

\bibitem[H\'{e}e]{He}
H\'{e}e,
{\it
Groupes de Chevalley et groupes classiques},
Seminar on finite groups, Vol. II, 1-54, Publ. Math. Univ. Paris VII, 17, Univ. Paris VII, Paris, 1984

\bibitem[J]{J}
Jantzen, J. C., Representations of Algebraic Groups, 2003
(American Math. Soc.)





\bibitem[KY07]{07} Kaneda M. and Ye, J., {\it Equivariant localization of $\bar D$-modules on the flag variety of the symplectic group of degree 4}, J. Alg. {\bf 309} (2007), 236--281

\bibitem[KY]{KY} Kaneda M. and Ye J.-C., {\it Some observatons on Karoubian complete strongly exceptional posets on the projective homogeneous varieties}, arXiv:0911.2568v1 [math.RT]



\bibitem[K09]{09} Kaneda M., {\it The structure of Humphreys-Verm modules for projective spaces}, J. Alg. {\bf 322} (2009), Pages 237-244 




\bibitem[K14]{14} Kaneda, M.,{\it Exceptional collections of sheaves on quadrics in positive characteristic}, S\~{a}o Paulo Journal of Mathematical Sciences {\bf 8} (2014), 117-156

\bibitem[K17]{17}
Kaneda, M,
{\it Another strongly exceptional collection of coherent sheaves on a Grassmannian}, 
Journal of Algebra 473 (2017) 352-373 




\bibitem[KaLa]{KaLa} Kashiwara M. and Lauritzen, N., {\it Local cohomology and $\cD$-affinity in positive characteristic}, C. R. Acad. Sci. Paris, Ser I {\bf 335}
(2002), 993-996.



\bibitem[Kat]{Kato} Kato S., {\it On the Kazhdan-Lusztig polynomials for affine Weyl groups}, Adv. Math. {\bf 55} (1995), 103-130


\bibitem[La]{La} Langer, A.,
{\it D-affinity
and Frobenius morphism on quadrics},
IMRN (2008), rnm 145



\bibitem[L80]{L80} Lusztig, G.,{\it Hecke algebras and Jantzen's generic decomposition patterns}, Adv. Math. {\bf 37} (1980),121-164



\bibitem[Or]{Or}
Orlov, D.O., {\it Derived categories of coherent
sheaves and equivalences between them}, Russian Math. Surveys. {\bf 58} (3) (2003), 511-591



\bibitem[St]{St} Steinberg, R., Lectures on Chevalley Groups, 2017 (AMS)


\bibitem[T]{T} Tanisaki, T., {\it Lie Algebras and Quantum Groups (in Japanese)
}, 2002
(共立出版)

\bibitem[Te]{Te} Testerman, D. M., Irreducible subgroups of exceptional algebraic groups, Memo {\bf 75}, No. 390, 1988 (AMS)























































































































\end{thebibliography}
\end{document}